\documentclass[11pt,a4paper,twoside]{article}
\usepackage{amsmath,amsfonts,amssymb,amsthm,bbm,latexsym,mathrsfs}
\usepackage{graphicx,color,epsfig,fancyhdr,dsfont}
\usepackage{enumerate}
\usepackage{hyperref}
\usepackage{indentfirst}
\usepackage[all]{hypcap}
\usepackage[affil-it]{authblk}
\allowdisplaybreaks
\title{Random Periodic Processes, Periodic Measures and Ergodicity}

\author[1]{Chunrong Feng}
\author[1]{Huaizhong Zhao}
\affil[1]{Department of Mathematical Sciences, Loughborough
University, LE11 3TU, UK}
\affil[ ]{C.Feng@lboro.ac.uk,  H.Zhao@lboro.ac.uk}
\date{}

\newtheorem{thm}{Theorem}[section]
\newtheorem{lem}[thm]{Lemma}
\newtheorem{cor}[thm]{Corollary}
\newtheorem{prop}[thm]{Proposition}
\newtheorem{defi}[thm]{Definition}
\newtheorem{rmk}[thm]{Remark}
\newtheorem{exmp}[thm]{Example}

\allowdisplaybreaks
\numberwithin{equation}{section}

\begin{document}

\maketitle

\begin{abstract}

Ergodicity of random dynamical systems with a periodic measure is obtained on a Polish space. In the Markovian case, the idea of Poincar\'e sections is introduced. It is proved that if the periodic measure is PS-ergodic, then it is ergodic. Moreover, if the infinitesimal generator of the Markov semigroup only has equally placed simple eigenvalues including $0$ on the imaginary axis, then the periodic measure is PS-ergodic and has positive minimum period. Conversely if the periodic measure with the positive minimum period is PS-mixing, then the infinitesimal generator only has equally placed simple eigenvalues (infinitely many) including $0$ on the imaginary axis. Moreover, under the spectral gap condition, PS-mixing of the periodic measure is proved. The ``equivalence" of random periodic processes and periodic measures is established. This is a new class of ergodic random processes. Random periodic paths of stochastic perturbation of the periodic motion of an ODE is obtained. 

\medskip

\medskip

\end{abstract}

\pagestyle{fancy}
\fancyhf{}
\fancyhead[LE,RO]{\thepage}
\fancyhead[LO]{\small{Random Periodic Processes, Periodic Measures and Ergodicity}}
\fancyhead[RE]{\small{C. R. Feng and  H. Z. Zhao}}


\section{Introduction}

Ergodicity is significant for the theory of random dynamical systems in describing their large time behaviour and irreducibility. 
However, important results 
have been proved 
only under the regime of stationary measures and 
stationary processes. They are not applicable to systems with periodicity.
In this paper we will break this restriction to provide an ergodic theory when ``periodicity" exists.
This scenario, regarded as random periodic, is defined in a very general situation of a separable Banach space, applicable for both
discrete random mappings and continuous time stochastic flows.



It is well known but still worth mentioning in this context  that
the notion of periodic paths has been a major concept in the theory of
dynamical systems since Poincar\'e's pioneering work (\cite{poincare}). 
Moreover, periodic phenomena exist in
many real world problems.
But, by nature, many real world systems are very often subject to
the influence of internal or external  randomness.
Periodicity, nonlinearity 
and randomness are present and interweave in many real world phenomena. Random periodicity is ubiquitous and 
can be found, for example, in daily temperature variations, economic and business cycles,
internet traffic volumes, activity of sunspots and transition between ice age and interglacial period.

But periodicity and randomness do not seem to match each other naturally,
so the first major task is to define the random periodicity
in a general setting. The study of random periodicity has attracted considerable 
interests in literature recently. 

Physicists have attempted to study random perturbations to periodic solutions for some time. They used 
first order linear approximations or asymptotic expansions in a small noise regime, e.g. see
\cite{knobloch}.  
The approach in \cite{knobloch} was to seek $Y(t+\tau,\omega)$ returning to a neighbourhood of $Y(t,\omega)$ for each noise realisation, where
$\tau>0$ is a fixed number. This suggests that almost surely 
each sample path is in a neighbourhood of its mean, which is not far from the original 
unperturbed periodic path. This reveals certain 
information about the ``periodicity" under small noise perturbations. However, in many situations, the sample path may not always stay in a small 
neighbourhood of its mean even when noise is small. One of the obstacles to
make more progress was the lack of a rigorous mathematical definition and 
appropriate mathematical tools. There were some scattering attempts in mathematics literature raising and discussing 
random periodic orbits for time-one mappings (\cite{klunger}).
Our work provides a systematic approach applicable for both time-one mappings (\cite{lian}) and flows.

New observation was made in \cite {zh-zheng} which says that for each fixed $t$, $\{Y(k\tau+t,\omega)\}_{k\in {\mathbb N}}$
should be a stationary path of the $\tau$-mesh discrete random dynamical system, where $\tau$ denotes the period. 
This then led to the rigorous definition of random periodic paths and a series of 
new results (\cite{feng-zhao-zhou},\cite{feng-zhao},\cite{lian},\cite{zh-zheng}).
An alternative way to understand random periodic behaviour  is
to study periodic measures which describe periodicity in the sense of distributions (\cite{Ha}).
There are a few works in the literature attempting to study statistical solutions of certain types of SDEs with periodic forceings;
motivated in the context of studying the climate change problem when the seasonal cycle is taken into considerations (\cite{majda-2}); 
some mean field models in chemical 
reactions (\cite{scheutzow})
and Ornstein-Uhlenbeck processes (\cite{chojn},\cite{rockner1}). However, it seems that the periodic measure was written in the form (\ref{zhaoh61}) for the first time in \cite{zhao2014}\footnote{This paper is now replaced by the 
current paper and will not be submitted for a publication.} (see also \cite{wu1},\cite{flz1}). 
Our formulation includes an entrance law and time periodicity of the measure function.
We note here that the time periodicity of the measure function in (\ref{zhaoh61}) was suggested by Has'minskii in \cite{Ha}, 
but the entrance law was missing in his formulation except for discrete transition semigroup at the integral multiples of the period.  It is important to note from our work that both the periodic condition and 
entrance law at all time are not redundant in the definition given in \cite{zhao2014}.


The concept of random periodic paths has led to more progress on investigations of 
various issues in stochastic dynamics and modelling real world problems.
They include an intriguing observation of random periodic paths in the stochastic 
Timmerman-Jin model of El Nino phenomenon arising in climate dynamics (\cite{chekroun});
bifurcations of stochastic reaction diffusion equations (\cite{wang});
periodic random attractors of stochastic lattice systems (\cite{blw}); stochastic resonance (\cite{clrs}); strange attractors of a particular hyperbolic 
random dynamical systems where infinite number of random periodic paths were found (\cite{hl});
anticipating random periodic solutions of SDEs (\cite{wu1}); numerical 
analysis of random periodic solutions and periodic measures of SDEs (\cite{flz1},\cite{flz3}).

For Markovian random dynamical systems, we introduce the idea of Poincar$\rm\acute e$ sections $\{L_s\}_{s\geq 0}$ 
with $L_{s+\tau}=L_s$ such that for any $x\in L_t$, $P(s,x,L_{s+t})=1,\ s, t\geq 0$.
Thus $P(k\tau, x, L_s)=1$ when $x\in L_s$. Note at integral multiples of the period $k\tau$, the discrete $\tau$-mesh random dynamical system 
$\{\Phi(k\tau)\}_{k\in {\mathbb N}}$ and 
its transition probability $\{P(k\tau)\}_{k\in {\mathbb N}}$, are in a 
stationary regime on each Poincar\'e section. We can apply  the Krylov-Bogoliubov procedure and the Chapman-Kolmogorov equation to find invariant measure $\rho_s$ 
with respect to $\{P(k\tau)\}_{k\in {\mathbb N}}$ on each Poincar\'e section $L_s$. These $\{\rho_s\}_{s\geq 0}$, form a periodic measure with respect to $\{P(t)\}_{t\geq 0}$. 
Moreover, if $\{P(k\tau)\}_{k\in {\mathbb N}}$ is irreducible on each $L_s$, then we can prove that the Poincar\'e sections are uniquely determined.

For a periodic measure $\{\rho_s\}_{s\in {\mathbb R}}$, its average over a period $\bar\rho={1\over \tau}\int _0^{\tau}\rho_sds$ is an invariant measure with respect to $\{P(t)\}_{t\geq 0}$. 
Thus we can construct a canonical dynamical system on a path space from the invariant measure, of which the ergodicity defines that of the invariant measure $\bar\rho$.
The periodic measure  $\{\rho_s\}_{s\in {\mathbb R}}$ is called ergodic if $\bar\rho$ is ergodic. 
 In a non-degenerate random periodic regime, 
the invariant measure $\bar \rho$ cannot weakly mixing and the transition probability $P(t,\cdot,\Gamma)$ does not converge. 
However, we will prove,
\begin{eqnarray*}
{\mathcal R}_N(\Gamma):=\int _{{\mathbb X}}|{1\over N}\sum_{k=0}^{N-1}\int _{0}^{\tau} (P(k\tau+s,x,\Gamma)-\rho_s(\Gamma))ds|\bar\rho(dx) \to 0,
\end{eqnarray*}
as $N\to \infty$ for each $\Gamma \in {\mathcal B}({\mathbb X})$ if and only if  the periodic measure is ergodic.

The concept of Poincar\'e sections is a key tool for establishing criteria for the convergence of ${\mathcal R}_N$. 
Observe that for each $s\in {\mathbb R}$, $\rho_s$ is an invariant measure of the $\tau$-mesh 
discrete Markovian semigroup $\{P(k\tau)\}_{k\in {\mathbb N}}$
on the Poincare section $L_s$. We call the periodic measure is PS-ergodic (PS-mixing) if 
for each $s$, the measure $\rho_s$ as an invariant measure of $\{P(k\tau)\}_{k\in {\mathbb N}}$ on $L_s$ is ergodic (mixing). 
We will prove that if the periodic measure is PS-ergodic, then it is ergodic.

We will 
classify between a real random periodic
regime and a degenerate stationary case.  
In the case of non-degenerate periodic measure with a minimum period $\tau>0$, 
 there is an angle variable which is not constant, unlike the stationary case. Thus the transformation operator has infinitely many eigenvalues on the unit circle.
In particular,  if the infinitesimal generator of the Markov semigroup only has equally placed simple eigenvalues including $0$ on the imaginary axis, then the periodic measure is PS-ergodic and has positive minimum period. Conversely if the periodic measure with the positive minimum period is PS-mixing, then the infinitesimal generator only has equally placed simple eigenvalues (infinitely many) including $0$ on the imaginary axis.  This is clearly distinguished from the 
mixing stationary case in which the Koopman-von Neumann Theorem says there is only one simple eigenvalue $0$ on the imaginary axis.

It is noted that the spectral structure of the Markov semigroup is more fruitful than that of the transformation operator on the path space.
In this context, it is worthy mentioning that in the case of the stationary regime, many results on
spectral gaps have been obtained, which give how far the rest of spectra of the generator are away from the $0$ eigenvalue (c.f. see \cite{chen},\cite{wangf} etc). 
Moreover, the spectral gap gives mixing property and convergence rate of transitional probability to the invariant measure. 
This fundamentally important result has brought many powerful analysis tools to the study of ergodicity 
and mixing of the invariant measure of stochastic systems.
In this paper, we prove if the semigroup has a spectral gap on each Poincar\'e section, then
the periodic measure is PS-mixing and for any $\Gamma \in {\mathcal B}({\mathbb R}^d)$, as $k\to \infty$,
\begin{eqnarray*}
\int _{{\mathbb X}}|{1\over \tau}\int _{k\tau}^{(k+1)\tau} P(s,x,\Gamma)ds-\bar\rho (\Gamma)|\bar\rho(dx)\leq {\rm e}^{-\delta k\tau}.
\end{eqnarray*}
This result, together with the result of eigenvalues on the imaginary axis, provides a clear 
analytic characterisation
of the PS-mixing property 
in terms of the spectra of the corresponding semigroup.
However, 
it is still open to see whether or not the spectral gap of differential operator implies the spectral gap of its semigroup on Poincar\'e sections.

Our ergodic theory give an innovative insight into the stochastic resonance and reveals a rigorous proof of the transition between ice age and interglacial period proposed in \cite {BenziStochRes} and the
partial differential equation for expected transition times (\cite{fzz},\cite{fzz2}).

Random periodic paths describe random periodicity in a pathwise manner, while periodic measure gives a  description
in terms of the law. 
They are not immediately equivalent, but both indispensable
concepts for understanding random periodicity, as stationary processes and invariant measures for the stationary regime.  
In this paper, we will prove that they can be "equivalent" in the following sense. 
First random periodic paths  give rise to
a periodic measure and conversely we are able to construct an enlarged probability space by adding trajectories of the random dynamical systems to be 
part of the new noise paths space, in which we can construct random periodic paths.
Moreover, one can prove that the law of the random periodic paths is the very periodic measure.

We would also like to point out that what we normally observe in the real world is only one realisation of a random periodic process, rather than a periodic measure.
However, random periodicity could be difficult to statistically test without appealing to the periodic distribution idea, especially when noise is large. On the other hand, to find a periodic measure from one realisation is an difficulty. 
To overcome this difficulty, we appeal to establish the law of large numbers and central limit theorem. We will publish these results in a different publication (\cite{flz2}).

\section{Random periodic paths and examples}\label{section 2.1}


%


In this part, we will study random periodicity of random dynamical systems of cocycles. 
This is necessary because on one hand random periodicity exists naturally for systems of cocycles. 
In this case, the integration of a periodic measure, if exists, over the time of one period is an invariant measure. Thus
its ergodicity makes sense as that of the average invariant measure.  On the other hand however, the above observation is not valid for 
stochastic periodic semi-flows. One cannot define an invariant measure from the integration of periodic measures thus ergodicity cannot be defined 
in the same way as above. But in the second part of this paper, we will use the idea of lifting stochastic periodic semi-flows to a cocycle on a cylinder, 
and periodic measure to that of the cocycle on the cylinder, of which the ergodicity can be studied. Thus the first fundamental 
task is to study the ergodicity of coycles.

Let $\mathbb{\mathbb{X}}$ be a Polish space and ${\mathcal B}(\mathbb{\mathbb{X}})$ be its Borel $\sigma$-algebra.
In this section, we
consider a measurable cocycle random dynamical system $\Phi$ on $(\mathbb{\mathbb{X}},{\mathcal B}(\mathbb{\mathbb{X}}))$ over a metric
dynamical system $(\Omega, {\mathcal F}, P, (\theta(t))_{t\in \mathbb{R}^+})$ with a one-sided time set $\mathbb{R}^+:=\{t\in \mathbb R:t\geq 0\}$, $\Phi: \mathbb{R}^+\times \Omega\times
\mathbb{\mathbb{X}}\to \mathbb{\mathbb{X}}$. It is $({\mathcal B}_{\mathbb{R^+}} \otimes {\mathcal{F}}\otimes {\mathcal{B}(\mathbb{\mathbb{X}})}, {\mathcal {B}(\mathbb{\mathbb{X}})})$-measurable and satisfies the cocycle property: $$\Phi(0,\omega)=id_{\mathbb{\mathbb{X}}} \ {\rm and } \
\Phi(t+s,\omega)=\Phi(t,\theta(s)\omega)\Phi(s,\omega), \ {\rm for \ all} \ s,t\in \mathbb{R}^+,
$$
for almost all $\omega\in \Omega$. The
map $\theta(t):\Omega\to \Omega$ is
$P$-measure preserving and measurably invertible. Therefore it can be extended to $\mathbb{R}^-$ as well by setting $\theta (t)=\theta (-t)^{-1}$ when $t\in \mathbb{R}^-$.
There is no need to require the map $\Phi (t,\omega): \mathbb{\mathbb{X}}\to \mathbb{\mathbb{X}}$ to be invertible, thus our work is applicable to both SDEs and SPDEs.
 
 First recall the definition of random periodic paths given in \cite{feng-zhao-zhou},\cite{feng-zhao}.
 
 \begin{defi}\label{zhao20179d}
A random periodic path of period $\tau$ of the random dynamical system $\Phi: \mathbb R^+\times \Omega\times \mathbb{\mathbb{X}}\to \mathbb{\mathbb{X}}$
 is an ${\mathcal F}$- measurable
map $Y:\mathbb{R}\times \Omega\to \mathbb{\mathbb{X}}$ such that for almost all $\omega\in \Omega$,
\begin{eqnarray}\label{prop6c}
\ \ \ \ \ \ \Phi(t,\theta(s)\omega) Y(s,\omega)=Y(t+s,\omega),\
Y(s+\tau,\omega)=Y(s, \theta(\tau) \omega),
\end{eqnarray}
for  any $t\in \mathbb{R}^+,s\in \mathbb{R}.$
It is called a random periodic path with the minimal period 
$\tau$ if $\tau>0$ is the smallest number such that (\ref{prop6c}) holds. 
It is a stationary path of $\Phi$ if $Y(s,\theta(-s)\omega)=Y(0,\omega)=:Y_0(\omega)$ for all $s\in \mathbb{R},\ \omega\in \Omega$, i.e. $Y_0:\Omega\to \mathbb{\mathbb{X}}$ is a stationary path
 if
for almost all $\omega\in \Omega$,
\begin{eqnarray*}
\Phi(t,\omega) Y_0(\omega)=Y_0(\theta (t)\omega), \ for\  any \ t\in \mathbb{R}^+.
\end{eqnarray*}
\end{defi}

The first part of the definition of the random periodic path suggests that a random periodic path $\{Y(s,\omega)\}_{s\in {\mathbb R}}$ is indeed a pathwise trajectory of the random dynamical system.
The second part of the definition says that it has some periodicity. But it is different from a periodic
path in the deterministic case, $Y(s+\tau,\omega)$
is not equal to $Y(s,\omega)$, but $Y(s,\theta (\tau)\omega)$. We call this random periodicity.
Starting at $Y(s,\omega)$, after a period $\tau$, trajectory does not return to $Y(s,\omega)$, but to $Y(s,\cdot)$ with different realisation $\theta ({\tau})\omega$. 
So it is neither completely
random, nor completely periodic, but a mixture of randomness and periodicity. In fact,  
the path $\{Y(s+\tau,\omega)\}_{s\in {\mathbb R}}$ repeats the path $\{Y(s,\theta (\tau)\omega)\}_{s\in {\mathbb R}}$, rather than
$\{Y(s,\omega)\}_{s\in {\mathbb R}}$ as in the deterministic case. This kind of random periodicity can be numerically 
checked as demonstrated in \cite{flz1}.

Let $\phi(s,\omega)=Y(s,\theta(-s)\omega), s\in \mathbb{R}$. It is easy to see that $\phi$ satisfies the definition in  \cite{zh-zheng}
\begin{eqnarray}\label{prop6d}
\ \ \ \ \ \  \Phi(t,\omega) \phi (s,\omega)=\phi (t+s,\theta (t)\omega),\ 
\phi (s+\tau,\omega)=\phi (s,  \omega), t\in \mathbb{R}^+, s\in \mathbb{R}.
\end{eqnarray}
Therefore $\phi$ is a periodic function and define
\begin{eqnarray}\label{prop6e}
L^{\omega}=\{(\phi (s,\omega): s\in [0,\tau)\}.
\end{eqnarray}
 It is easy to see from the first formula in (\ref{prop6d}) that $L^{\omega}$ is an invariant set, i.e.
 $$\Phi(t,\omega)L^{\omega}=L^{\theta (t)\omega},$$
 for any  $t\in \mathbb{R}^+$.  But needless to say that random periodic solution
 gives more detailed information about the dynamics of the random dynamical system than a general invariant set.
Unlike the periodic solution of deterministic dynamical systems,
the random dynamical system does not follow the closed curve, but move from one closed curved to another
when time evolves. This is fundamentally different from the deterministic case, which makes it hard to study. However, this natural definition in random case makes it possible to gain new
understanding of random phenomena with some periodic nature, where strict deterministic periodicity is not applicable.

The above definition is given for the continuous time case only. All the results are given in this setting as well. 
They all apply to the case when the time is discrete, i.e. when $\mathbb{R}$ is replaced by 
 $\mathbb{Z}=\{\cdots,-2,-1,0,1,2,\cdots\}$ and $\mathbb{R}^+$ by $\mathbb{N}=\{0,1,2,\cdots\}$.
 
It is not the purpose of this paper to discuss the existence of the random periodic path. 
In this paper, we only give one example of random dynamical systems that has a random periodic path.

%

The problem of a random perturbation to periodic motions is of great interests to both mathematicians and physicists.
If an ordinary differential equation (ODE) has a periodic path, does a stochastic differential equation with the coefficients of the ODE
as its drift possess a random periodic path? This can be regarded as stochastic perturbations of periodic motion of the dynamical system generated by the ODE. 
If the noise is nondegenerate (strictly elliptic), we can see that random periodic solution is 
synchronised to a stationary solution.

Zhao-Zheng (2009) provided a first example of stochastic differential equation with a random periodic path. This is SDE (\ref{zhao300x}) with $W_1=W_2$
instead of two independent Brownian motions. In this case, the random periodic path was written explicitly in Zhao-Zheng (2009). 
But when $W_1$ and $W_2$ are independent Brownian motions, SDE (\ref{zhao300x}) still has a random periodic path with a positive 
minimum period, but its proof is much more involved. Note the noise in  (\ref{zhao300x}) is degenerate. 

\begin{exmp}\label{zhao300v}
Consider the following stochastic differential equation on ${\mathbb R}^2$
\begin{eqnarray}\label{zhao300x}
\left\{
\begin{array}{cc}
dx_1=&[-x_2+x_1(1-x_1^2-x_2^2)]dt+x_1dW_1(t),\\
dx_2=&[x_1+x_2(1-x_1^2-x_2^2)]dt+x_2dW_2(t).
\end{array}
\right.
\end{eqnarray}
Here $W_1(t)$ and $W_2(t)$ are two independent one-dimensional two-sided Brownian motions on the probability space
$(\Omega, {\mathcal F}, P)$ with   $(W_1(0),W_2(0))^T=(0,0)^T$. Denote $W(t)=(W_1(t),W_2(t))^T$. 
Set ${\mathcal F}_s^t=\sigma(W(u)-W(v): s\leq v\leq u\leq t)$, ${\mathcal F}_{-\infty}^t=V_{s\leq t}{\mathcal F}_s^t$
and $\theta: {\mathbb R}\times \Omega\to \Omega$ the measure preserving metric dynamical system  
given by 
$$(\theta_s\omega)(t)=W(t+s)-W(s), s,t\in {\mathbb R}.
$$
Denote by $\Phi=(\Phi_1, \Phi_2):[0,\infty)\times {\mathbb R}^2\times\Omega\to {\mathbb R}^2$
the cocycle generated by the solutions of (\ref{zhao300x}).

It is well known that the noiseless system 
\begin{eqnarray}\label{13Jan7}
\left\{
\begin{array}{cc}
{dx_1\over dt}=&-x_2+x_1(1-x_1^2-x_2^2),\\
{dx_2\over dt}=&x_1+x_2(1-x_1^2-x_2^2),
\end{array}
\right.
\end{eqnarray}
has a periodic solution $(x_1(t),x_2(t))=(\cos t, \sin t)$. In the following 
proposition, we will study the existence of random periodic path which can be regarded as a random perturbation of the periodic motion of the deterministic 
dynamical system generated by (\ref{13Jan7}). 

\end{exmp}

\begin{prop}
Equation (\ref{zhao300x}) has a unique random periodic solution $x^*(t)=(x^*_1(t),x^*_2(t))\ne (0,0)$
with a positive minimum period satisfying for a.s. $\omega\in \Omega$,
\begin{eqnarray}
x^*(t+\pi,\omega)&=&-x^*(t,\theta _{\pi}\omega),\label{zhao310b}\\
x^*(t+2\pi,\omega)&=&x^*(t,\theta _{2\pi}\omega).\label{zhao310c}
\end{eqnarray}
\end{prop}
\begin{proof}
Let us use the polar coordinates 
by letting $x_1=r\cos\varphi, x_2=r\sin \varphi$. Then we have
\begin{eqnarray}\label{zhao300e}
\hskip30pt
\left\{
\begin{array}{ll}
dr=&r(1-r^2+{1\over 4}\sin^2(2\varphi))dt+r\cos ^2\varphi \ dW_1(t)+r\sin^2\varphi \ dW_2(t),\\
d\varphi=&dt+{1\over 4}\sin (4\varphi)dt+{1\over 2}\sin(2\varphi)\ d(W_2(t)-W_1(t)).
\end{array}
\right.
\end{eqnarray}
This generates a stochastic flow $(r,\varphi): [0,\infty)\times ({\mathbb R}^+\times {\mathbb R})\times\Omega\to ({\mathbb R}^+\times {\mathbb R})$.
Let us first look at the angle equation. Note that the coefficients $b(\varphi)=1+{1\over 4} \sin (4\varphi)$ and $\sigma (\varphi)={1\over 2} \sin(2\varphi)$ are periodic functions of period ${\pi\over 2}$ and ${\pi}$ respectively. Thus we can consider the equation as an SDE on a circle of radius $1\over 2$ i.e. we can consider $\tilde\varphi_t=\varphi_t\mod \pi$, then $\tilde \varphi_t$ is a random dynamical system cocycle on the circle $S_{1\over 2}$. By the fact that $S_{1\over 2}$ is compact, so there is an invariant measure $\rho_{\tilde \varphi}$ for $\tilde \varphi$. Therefore by Birkhoff ergodic theorem, we have as $T\to \infty$, 
$${1\over T} \int_0^T \sin^2(2\varphi_t)dt={1\over T} \int_0^T \sin^2(2\tilde\varphi_t)dt\to \int_{S_{1\over 2}} \sin^2(2x)\rho_{\tilde \varphi}(dx).$$
When $\tilde \varphi=0$ or ${\pi\over 2}$, $d\varphi_t=d\tilde\varphi_t=dt$, then it is obvious that $\rho_{\tilde \varphi}$ 
cannot be supported at $\{0, {\pi\over 2}\}$. Thus 
\begin{eqnarray}\label{zhao310d}
\beta:=\int_{S_{1\over 2}}\sin^2(2x)\rho_{\tilde \varphi}(dx)>0, \ \ a.s.
\end{eqnarray}

Now we consider  $\psi_{t}=\varphi_t-t$, then 
\begin{eqnarray}\label{eqn1.10}
d_t \psi_{t}={1\over 4} \sin (4\psi_{t}+4t)dt+{1\over 2}\sin(2\psi_{t}+2t)d(W_t^2-W_t^1).
\end{eqnarray}
Denote by $\psi_t(\alpha)$ as the solution of (\ref{eqn1.10}) 
with initial condition $\psi_{0}=\alpha$. Note $\psi_t$ satisfies a stochastic differential equation with coefficients periodic in time with period $\pi$.
Inspired by Carvehille-Chappell-Elworthy \cite{cce} (see also Rogers-Williams \cite{rw}), we consider the gradient flow on the circle $S_{1\over 2}$ and its Lyapunov exponent. 
Define $\Psi_{t}(\alpha)=\nabla _{\alpha} \psi_{t}(\alpha).$ Then it is easy to see that 
\begin{eqnarray*}
d_t\Psi_{t}=\cos(4\psi_{t}+4t)\Psi_{t}dt+\sqrt {2} \cos(2\psi_{t}+2t)\Psi_{t} d({{W_t^2-W_t^1}\over {\sqrt 2}}).
\end{eqnarray*}
This is a linear stochastic differential equation for $\Psi_t$. Note that ${W_t^2-W_t^1}\over {\sqrt 2}$ is equivalent to a standard one-dimensional Brownian motion, so
by It$\hat {\rm o}$'s formula, we can solve
\begin{eqnarray*}
\Psi_{t}(\alpha)&=&\exp \{\int_0^t \cos(4\psi_{r}+4r)dr-\int_s^t \cos^2(2\psi_{r}+2r)dr\\
&&\hskip1cm +\int_0^t \cos(2\psi_{r}+2r)d(W_r^2-W_r^1)\}\\
&=&\exp \{-\int_0^t \sin^2(2\psi_{r}+2r)dr+\int_0^t \cos(2\psi_{r}+2r)d(W_r^2-W_r^1)\}.
\end{eqnarray*}
Thus the Lyapunov exponent is computed as follows by using (\ref{zhao310d}),
\begin{eqnarray*}
\lambda&=&\lim_{t\to \infty}{1\over {t}}\log |\Psi_{t}(\alpha)|\\
&=&-\lim_{t\to \infty}{1\over {t}}\int_s^t \sin^2(2\varphi_r)dr\\
&=&-\beta<0, \ \ a.s.
\end{eqnarray*}
Then there exists a tempered random variable $C(\omega)>0$ such that  for a.s. $\omega\in \Omega$
$$|\psi_{s+k\pi}(\theta _{-k\pi}\omega, \alpha)-\psi_{s+k\pi}(\theta _{-k\pi}\omega, \alpha')|\leq C(\omega)
e^{-\beta k\pi}.$$ In particular, for a.s. $\omega$, for $k<l$, 
$$|\psi_{s+k\pi}(\theta _{-k\pi}\omega, \alpha)-\psi_{s+l\pi}(\theta _{-l\pi}\omega, \alpha)|\leq C(\omega)e^{-\beta k\pi}.$$ 
Thus $\{\psi_{s+k\pi}(\theta _{-k\pi}\omega, \alpha)\}_k$ is a Cauchy sequence and therefore it has a limit, denoted by $\psi_s^*(\omega)$. The limit does not depend on $\alpha$. Note, for $t\geq 0$, 
\begin{eqnarray*}
\psi_{t+s+k\pi}(\theta _{-k\pi}\omega, \alpha)=\psi_{t}(\theta _s\omega)\circ \psi_{s+k\pi}(\theta _{-k\pi}\omega, \alpha)\to \psi_{t}(\theta _s\omega)\psi_s^*(\omega).
\end{eqnarray*}
But for a.s. $\omega\in \Omega$
\begin{eqnarray*}
\psi_{t+s+k\pi}(\theta _{-k\pi},\alpha)\to \psi_{t+s}^*(\omega).
\end{eqnarray*}
Thus for a.s. $\omega\in \Omega$, $t\geq 0$,
$$\psi_{t}(\theta _s\omega)\psi_s^*(\omega)=\psi_{t+s}^*(\omega).$$
Moreover, for a.s. $\omega\in \Omega$
$$\psi_{s+\pi+k\pi}(\theta _{-k\pi}\omega,\alpha)=\psi_{s+(k+1)\pi}(\theta _{-(k+1)\pi}\theta _{\pi}\omega,\alpha)\to \psi^*_{s}(\theta_\pi\omega)$$
and for a.s. $\omega\in \Omega$
$$\psi_{s+\pi+k\pi}(\theta _{-k\pi}\omega, \alpha)\to \psi^*_{s+\pi}(\omega).$$
Thus for a.s. $\omega\in \Omega$
$$\psi^*_{s+\pi}(\omega)=\psi^*_{s}(\theta_\pi\omega).$$
This means that SDE (\ref{eqn1.10}) has a random periodic solution with period $\pi$. Set 
\begin{eqnarray}
\label{zhao310e}
\varphi^*_{s}(\omega)=s+\psi^*_{s}(\omega), \ {\rm for\ any\ s\in \mathbb R}.
\end{eqnarray}
Then it is easy to see that 
\begin{eqnarray}
\varphi^*_{s+\pi}(\theta_{-\pi}\omega)=\pi+\varphi^*_{s}(\omega)\label{1.11}\\
\varphi^*_{s+2\pi}(\theta_{-2\pi}\omega)=2\pi+\varphi^*_{s}(\omega).\label{1.12}
\end{eqnarray}
Moreover, let $\varphi_s(\cdot)$ denote the stochastic flow generated by the second equation of SDE (\ref{zhao300e}). Then for $t\geq 0$, 
\begin{eqnarray}\label{1.13}
\varphi_t(\theta_s\omega)\varphi^*_s(\omega)=\varphi^*_{t+s}(\omega).
\end{eqnarray}
Consider SDE (\ref{zhao300e}), 
the radius and angle coordinates together generate a cocycle satisfying for $s, t\geq 0$, 
$$(r_t(\theta_s\omega), \varphi_t(\theta_s\omega))\circ (r_s(\omega), \varphi_s(\omega))=(r_{t+s}(\omega), \varphi_{t+s}(\omega)).$$
On the other hand, inspired by Arnold \cite{ar},
let $\xi={1\over {r^2}}$, then 
\begin{eqnarray}\label{zhao310a}
d\xi&=&\xi(-3+4(\cos^4\varphi+\sin^4\varphi)) dt+2dt\nonumber 
\\
&&+\xi(-2\cos^2\varphi dW_1(t)-2\sin^2\varphi dW_2(t)).
\end{eqnarray}
Denote by $\varphi_{t,s}(\varphi),\xi_{t,s}(\xi,\varphi)$ the solution of the second equation in SDE (\ref{zhao300e}) and the solution of (\ref{zhao310a})
respectively for $t\ge s$ with $\varphi_{s,s}(\varphi)=\varphi$ and $\xi(s,s,\xi,\varphi)=\xi$.
Then given $(\xi_s, \varphi_s)$ being measurable with respect to ${\mathcal F}_{-\infty}^s$,  one can solve $\xi$ easily as follows, for $t\geq s$,
\begin{eqnarray*}
&&\xi(t, s, \xi_s,\varphi_s)\\
&=&\xi_s \exp \{-3(t-s)+2\int_s^t (\cos ^4 (\varphi_{v,s}(\varphi_s))+\sin ^4 (\varphi_{v,s}(\varphi_s)))dv\\
&&\hskip1.2cm-2\int_s^t \cos ^2 (\varphi_{v,s}(\varphi_s))dW_1(v)-2\int_s^t \sin ^2 (\varphi_{v,s}(\varphi_s))dW_2(v)\}\\
&&+2 \int_s^t \exp \{-3(t-u)+2\int_u^t (\cos ^4 (\varphi_{v,s}(\varphi_s))+\sin ^4 (\varphi_{v,s}(\varphi_s)))dv\\
&&\hskip 2cm -2\int_u^t \cos ^2 (\varphi_{v,s}(\varphi_s))dW_1(v)-2\int_u^t \sin ^2 (\varphi_{v,s}(\varphi_s))dW_2(v)\}\ du.
\end{eqnarray*}
Given $(r_s,\varphi_s)$, consider $(\xi_s, \varphi_s)=({1\over {r_s^2}}, \varphi_s)$.
It follows that for $t\geq s$, 
\begin{eqnarray*}
&&r(t,s, r_s, \varphi_s)\\
&=& (\xi(t,s,\xi_s,\varphi_s))^{-{1\over 2}}\\
&=&\left[ \exp\{3(t-s)-2\int_s^t (\cos ^4 (\varphi_{v,s}(\varphi_s))+\sin ^4 (\varphi_{v,s}(\varphi_s)))dv\right .\\
&&\left .\hskip1cm +2\int_s^t (\cos ^2 (\varphi_{v,s}
(\varphi_s))dW_1(v)+ \sin ^2 (\varphi_{v,s}(\varphi_s))dW_2(v))\}\right ]^{1\over 2}
\\
&&\cdot \left [{1\over {r_s^2}}+2\int_s^t \exp\{3(u-s)-2\int_s^u (\cos ^4 (\varphi_{v,s}(\varphi_s))+\sin ^4 (\varphi_{v,s}
(\varphi_s)))dv\right .\\
&&
\left .\hskip1cm +2\int_s^u (\cos ^2 (\varphi_{v,s}(\varphi_s))dW_1(v)+ \sin ^2 (\varphi_{v,s}(\varphi_s))dW_2(v))\}du\right]^{-{1\over 2}},
\end{eqnarray*}
defines the solution of the first equation of (\ref {zhao300e}) for any $t\geq s$ wth initial condition $r_s$ at the time $t=s$.
Recall $\varphi^*$ given in (\ref{zhao310e}). Define, for any $s\in \mathbb R$, 
\begin{eqnarray*}
r_s^*
&:=&\left[2\int_{-\infty}^s\exp\{-3(s-u)+2\int_u^s (\cos ^4 \varphi_v^*+\sin ^4 \varphi_v^*)dv\right .\\
&& \left .\hskip1cm -2\int_u^s (\cos ^2 \varphi_v^*dW_1(v)+ \sin ^2 \varphi_v^*dW_2(v))\}du\right]^{-{1\over 2}}.
\end{eqnarray*}
Then for any $t\geq s$,
\begin{eqnarray}
&&r(t, s, r_s^*, \varphi_s^*,\omega)\nonumber\\
&=& \Big[\int_{-\infty}^s\exp\{-3(t-u)+2\int_u^t (\cos ^4 \varphi_v^*+\sin ^4 \varphi_v^*)dv
\nonumber \\
&&\hskip2cm 
-2\int_u^t (\cos ^2 \varphi_v^*dW_1(v)+\sin ^2 \varphi_v^*dW_2(v))\}du\nonumber\\
&&+\int_s^t\exp\{-3(t-u)+2\int_u^t (\cos ^4 \varphi_v^*+\sin ^4 \varphi_v^*)dv
\nonumber \\
&&\hskip2cm 
-2\int_u^t (\cos ^2 \varphi_v^*dW_1(v)+ \sin ^2 \varphi_v^*dW_2(v))\}du\Big]^{-{1\over 2}}\label{1.14}\\
&=&\Big[\int_{-\infty}^t\exp\{-3(t-u)+2\int_u^t (\cos ^4 \varphi_v^*+\sin ^4 \varphi_v^*)dv
\nonumber \\
&&\hskip2cm 
-2\int_u^t (\cos ^2 \varphi_v^*dW_1(v)+ \sin ^2 \varphi_v^*dW_2(v))\}du\Big]^{-{1\over 2}}\nonumber\\
&=&r_t^*(\omega). \nonumber
\end{eqnarray}
Moreover, let $\tilde W_i(v)=\theta_\pi W_i(v)=W_i({v+\pi})-W_i(v), i=1, 2$, then by the change of variables and (\ref{1.11}),
\begin{eqnarray}
\label{1.15}
&&
r^*_{s+\pi}(\omega)\\
&=&\Big[\int_{-\infty}^{s+\pi}\exp\{-3({s+\pi}-u)+2\int_u^{s+\pi} (\cos ^4 \varphi_v^*+\sin ^4 \varphi_v^*)dv\nonumber \\
&&\hskip2cm -2\int_u^{s+\pi} (\cos ^2 \varphi_v^*dW_1(v)+ \sin ^2 \varphi_v^*dW_2(v))\}du\Big]^{-{1\over 2}}\nonumber\\
&=&\Big[\int_{-\infty}^{s+\pi}\exp\{-3({s+\pi}-u)+2\int_{u-\pi}^{s} (\cos ^4 \varphi^*_{v}(\theta_{\pi}\omega)+\sin ^4 \varphi^*_{v}(\theta_{\pi}\omega))dv
\nonumber \\
&&\hskip2cm -2\theta_\pi\circ \int_{u-\pi}^{s} (\cos ^2 \varphi_{v}^*d W_1(v)+\sin ^2 \varphi_{v}^*d W_2(v))\}du\Big]^{-{1\over 2}}\nonumber\\
&=&\Big[\int_{-\infty}^{s}\exp\{-3({s}-u)+2\int_{u}^{s} (\cos ^4 \varphi^*_{v}(\theta_{\pi}\omega)+\sin ^4 \varphi^*_{v}(\theta_{\pi}\omega))dv\nonumber \\
&&\hskip2cm -2\theta_\pi\circ \int_{u}^{s} (\cos ^2 \varphi_{v}^*d W_1(v)+\sin ^2 \varphi_{v}^*d W_2(v))\}du\Big]^{-{1\over 2}}\nonumber\\
&=& r^*_s(\theta_{\pi}\omega).\nonumber
\end{eqnarray}
That is to say that $r^*$ is random periodic with period $\pi$. 
Let $$(x_1^*(s), x_2^*(s))=(r^*_s \cos \varphi_s^*, r^*_s\sin (\varphi^*_s)).$$
Then from (\ref{1.11}) and (\ref{1.15}), we know that 
\begin{eqnarray*}
(x_1^*(s+\pi,\omega), x_2^*(s+\pi,\omega))&=&(r^*_{s+\pi}\cos \varphi_{s+\pi}^*(\omega), r^*_{s+\pi}
\sin \varphi_{s+\pi}^*(\omega))\\
&=& -(r^*_s(\theta_\pi\omega)\cos (\varphi_{s}^*(\theta_\pi\omega)), r^*_s(\theta_\pi\omega)\sin (\varphi_{s}^*(\theta_\pi\omega)))\\
&=&-(x_1^*(s,\theta_\pi\omega), x_2^*(s,\theta_\pi\omega)),
\end{eqnarray*}
i.e. (\ref{zhao310b}) holds. Similarly, by (\ref{1.12}) and (\ref{1.15}), we can prove that (\ref{zhao310c}) also holds.
Now by (\ref{1.13}) and (\ref{1.14}), we know that for $t\geq 0$, 
\begin{eqnarray*}
&&\Phi(t, \theta_s\omega)(x^*_1(s,\omega), x^*_2(s,\omega))\\
&=& (\Phi_1(t, \theta_s\omega), \Phi_2(t,\theta_s\omega))(x^*_1(s,\omega), x^*_2(s,\omega))\\
&=&(r^*_{t+s}(\omega)\cos(\varphi^*_{t+s}(\omega)), r^*_{t+s}(\omega)\sin(\varphi^*_{t+s}(\omega)))\\
&=&(x^*_1(t+s,\omega), x^*_2(t+s,\omega)).
\end{eqnarray*}
That is to say we have a random periodic solution $(x^*_1(s,\omega), x^*_2(s,\omega))\neq (0,0)$, $s\in \mathbb R$ with periodic $2\pi$. It is clear from (\ref{zhao310b}) that 
the minimum period is strictly positive.
\end{proof}

\begin{rmk} 
%

(i). We have done some numerical simulations to equation (\ref{zhao300x}). To explain the numerical simulations, 
note (\ref{prop6c}) (or (\ref{zhao310c})) implies 
$$x^*(s-2\pi,\omega)=x^*(s,\theta_{-2\pi}\omega) {\rm \ for \ all }\ s\in R.
$$
 This means the paths $x^*(\cdot,\omega)$ and $x^*(\cdot,\theta_{-2\pi}\omega)$ are identical if we shift each coordinate of $x^*(\cdot,\theta_{-2\pi}\omega)$
to the left by $\tau$. 
By the same reason, if $x^*(\cdot, \omega)$ is a stationary path,  then for any $t$, $x^*(\cdot,\omega)$ and $x^*(\cdot,\theta_{-t}\omega)$ are identical if we shift each coordinate of $x^*(\cdot,\theta_{-t}\omega)$
to the left (when $t>0$) or the right (when $t<0$) by $|t|$, since a stationary path is when (\ref{prop6c}) holds for $\tau$ being any real number. 
  The numerical simulations demonstrated in Figure 1 describe that $x^*_1(\cdot,\omega)$ and $x^*_1(\cdot,\theta_{-2\pi}\omega)$ are identical up to a shift, while $x^*_1(\cdot,\theta_{-\pi}\omega)$ is not identical to them. But our simulations suggest $x^*_1(\cdot-\pi,\omega)=-x^*_1(\cdot,\theta_{-\pi}\omega)$, which is exactly what we proved in Proposition 1.1.4.
  This provides numerical evidence that $x^*$ is not a stationary path.
   
  (ii). It is obvious from (\ref{zhao310b}) and (\ref{zhao310c}) that $(x_1^*,x_2^*)$ is the random periodic path with a positive minimum period. 
 It is also easy to draw the conclusion that ${2\pi\over n}$, $n$ being even members, cannot be the minimum period. The minimum period has to be of the form 
 $2\pi\over n$ with $n$ being an odd number. But it is not clear 
whether or not $2\pi$ is indeed its minimum period. 
   
  (iii). To compare the situation with a stationary solution case, we consider a similar perburbed equation with additive noise:
  \begin{eqnarray}\label{zhao20177e}
\left\{
\begin{array}{cc}
dx_1=&[-x_2+x_1(1-x_1^2-x_2^2)]dt+dW_1(t),\\
dx_2=&[x_1+x_2(1-x_1^2-x_2^2)]dt+dW_2(t).
\end{array}
\right.
\end{eqnarray}
 This  equation has a stationary path. Indeed numerical simulations demonstrate that $x^*_1(\cdot,\omega)$, $x^*_1(\cdot,\theta_{-\pi}\omega)$ and $x^*_1(\cdot,\theta_{-2\pi}\omega)$ are identical up to a shift (Figure 2). We have done simulations of pull-back of some other values of time as well. Though not presented here for the interests of space, they are all identical up to a shift.
  

\begin{figure}
\begin{minipage}{\textwidth}
\centering 
\includegraphics[scale=0.20]{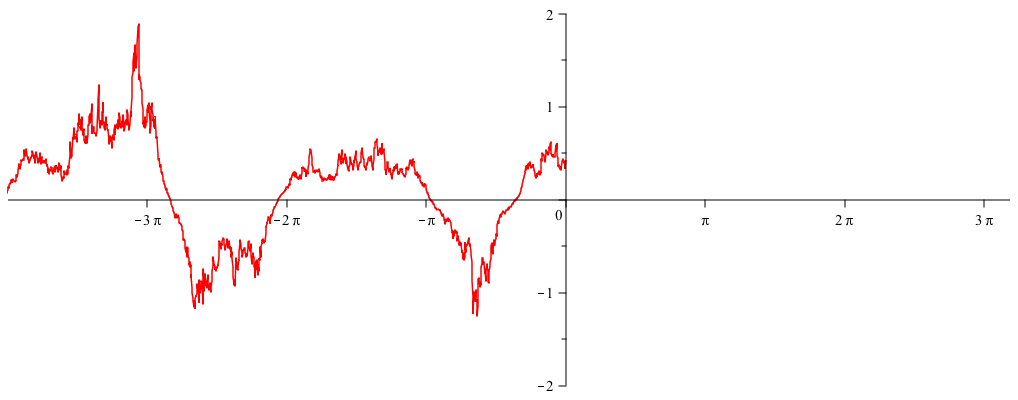}\\
\includegraphics[scale=0.20]{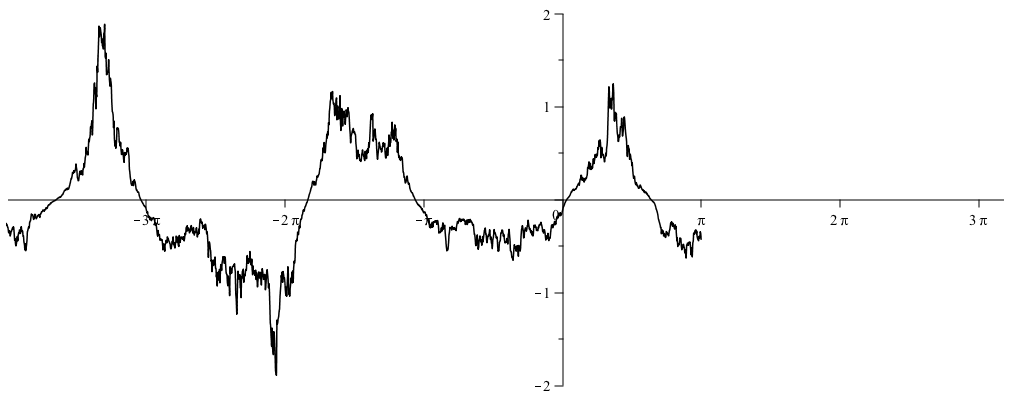}\\
\includegraphics[scale=0.20]{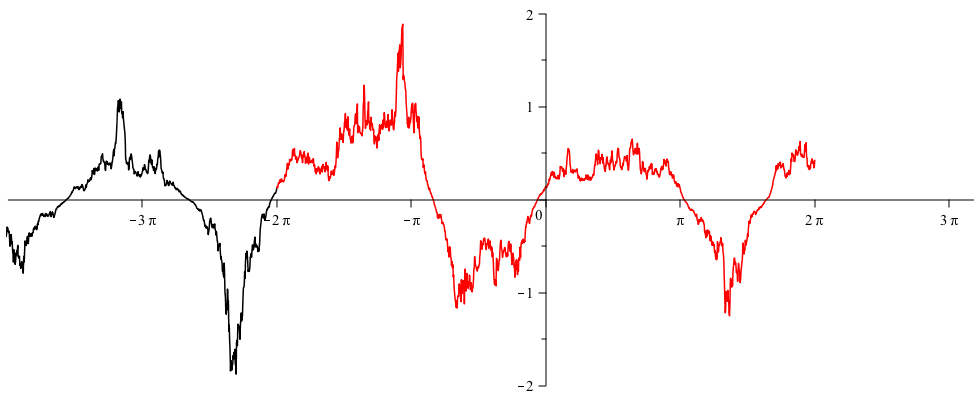}
\end{minipage} 
\caption{
Multiplicative noise-random periodic solutions: 
from the left to right, path of first coordinate with one realisation $\omega$, its pullbacks $\theta _{-\pi}\omega$ and $\theta_{-2\pi}\omega$ respectively. Red paths are identical up to a shift.
}
\label{graph of example 1 samples}
\end{figure}

\begin{figure}
\begin{minipage}{\textwidth}
\centering 
\includegraphics[scale=0.20]{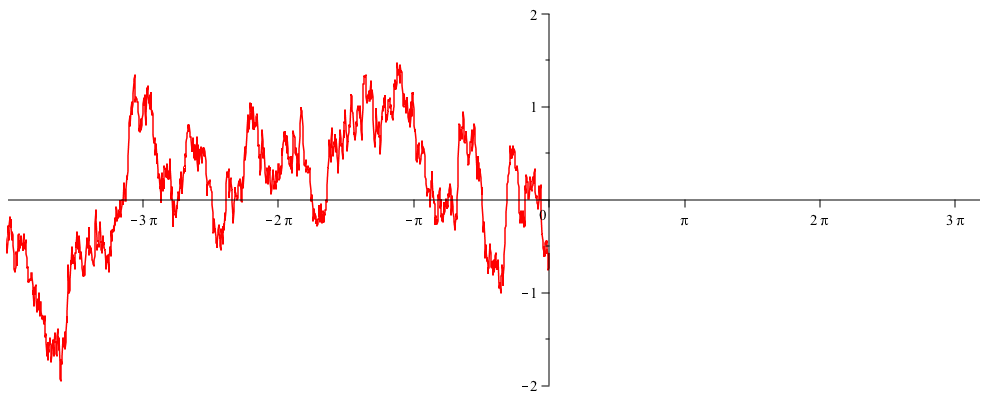}\\
\includegraphics[scale=0.20]{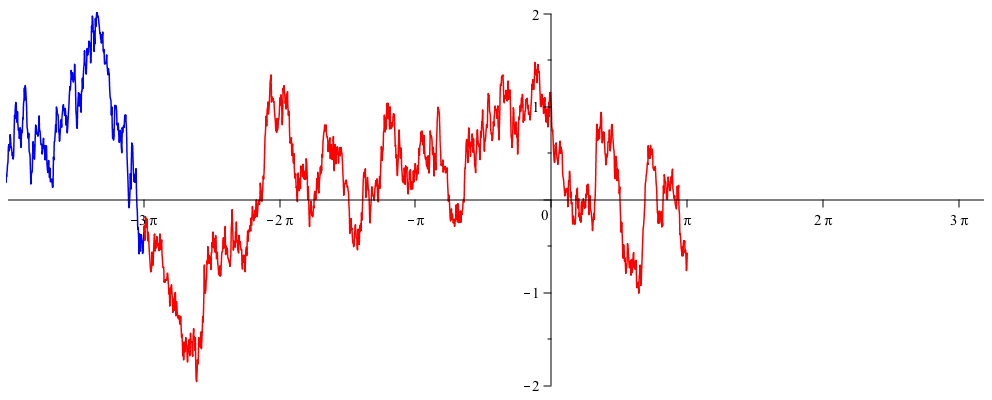}\\
\includegraphics[scale=0.20]{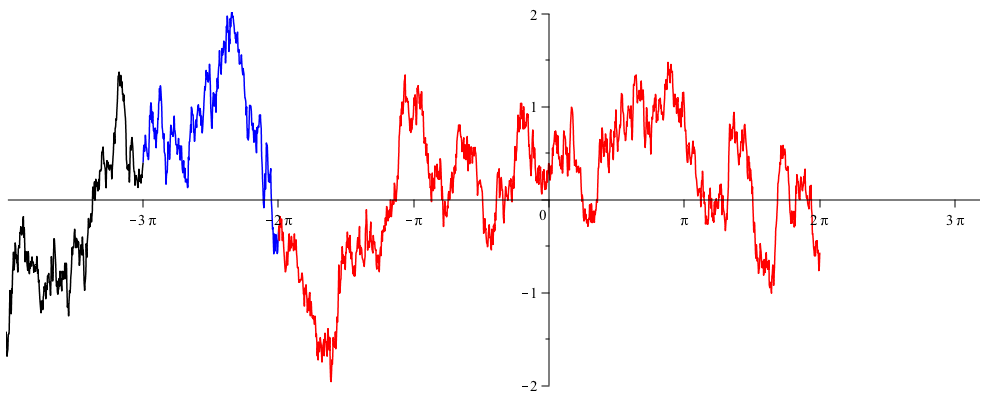}
\end{minipage} 
\caption{
Additive noise-stationary solutions: 
from the left to right, path of first coordinate with one realisation $\omega$, its pullbacks $\theta _{-\pi}\omega$ and $\theta_{-2\pi}\omega$ respectively. Red paths are identical up to a shift, so do blue paths.}
\label{graph of example 1 samples}
\end{figure}
\vskip5pt
\end{rmk}

\section{Periodic measures}

We start our investigation with proving a simple but important result that under the assumption of the existence of random periodic paths, the random Dirac measure with the support on sections of the random periodic curve $L^{\omega}$ is the periodic measure and its time average is an invariant measure. To make this clear, we consider a standard product measurable space $(\bar\Omega, \bar{\mathcal {F}})=(\Omega\times \mathbb{X}, {\mathcal {F}}\otimes {\mathcal {B}}(\mathbb{\mathbb{X}}))$ and the skew-product of the metric dynamical system $(\Omega,{\mathcal F},P,(\theta(t))_{t\in \mathbb{R}})$ and the cocycle $\Phi(t,\omega)$ on $\mathbb{\mathbb{X}}$,
$\bar {\Theta}(t):  \bar\Omega\to \bar\Omega$,
\begin{eqnarray}
\bar\Theta(t)(\bar\omega)=(\theta(t)\omega, \Phi(t,\omega)x), \
  {\rm where}\ \bar\omega=(\omega, x),\ t\in \mathbb{R}^+.
  \end{eqnarray}
Recall
\begin{eqnarray*}{\mathcal P}_P(\Omega\times \mathbb{\mathbb{X}}):=\{&\mu:&probability \ measure\  on\  (\Omega\times \mathbb{\mathbb{X}}, {\mathcal F}\otimes {\mathcal B}(\mathbb{\mathbb{X}}))\\
&&\ \ \ \ \ \ \ \ \ \ \ \ \ \ \ \  with\ marginal\  P \ on \ (\Omega, {\mathcal F})\}
\end{eqnarray*}
 and $${\mathcal P}(\mathbb{\mathbb{X}}):=\{\rho:  probability \ measure\  on\  (\mathbb{\mathbb{X}}, {\mathcal B}(\mathbb{\mathbb{X}}))\}.
 $$
\begin{defi}\label{zhaoh66}
A map $\mu: \mathbb R\to {\mathcal P}_P(\Omega\times \mathbb{\mathbb{X}})$ is called a periodic probability measure of period $\tau$ on $(\Omega\times \mathbb{\mathbb{X}}, {\mathcal {F}}\otimes {\mathcal {B}(\mathbb{\mathbb{X}})})$ for the random dynamical system $\Phi$  if
\begin{eqnarray} \label{1.3}
\mu_{\tau+s}=\mu_s \ and\ \bar\Theta(t)\mu_s=\mu_{t+s}, \ for\  any \ t\in \mathbb{R}^{+},s\in \mathbb{R}.
\end{eqnarray}
It is called a periodic measure with minimum period $\tau>0$ if $\tau$ is the smallest number such that (\ref{1.3}) holds.
It is an invariant measure if it also satisfies $\mu_s=\mu_0$ for any $s\in \mathbb{R}$ i.e.  $\mu_0$ is an invariant measure
of $\Phi$ if
$\mu_0\in {\mathcal P}_P(\Omega\times \mathbb{\mathbb{X}})$ and
\begin{eqnarray}
\bar\Theta(t)\mu_0=\mu_{0}, \ for\  any \ t\in \mathbb{R}^{+}.
\end{eqnarray}
\end{defi}
\begin{thm}\label{zhao13}
If a random dynamical system $\Phi: \mathbb{R}^+\times \Omega \times \mathbb{\mathbb{X}}\to \mathbb{\mathbb{X}}$ has a random periodic path $Y: \mathbb{R}\times \Omega \to \mathbb{\mathbb{X}}$, it has a periodic measure on $(\Omega\times \mathbb{\mathbb{X}}, {\mathcal {F}}\otimes {\mathcal{B}}(\mathbb{\mathbb{X}}))$ $\mu :\mathbb R\to {\mathcal P}_P(\Omega\times \mathbb{\mathbb{X}})$ given by
\begin{eqnarray}\label{1.4}
\mu_s(A)
=\int_\Omega \delta_{Y(s,\omega)}(A_{\theta(s)\omega})P(d\omega),
\end{eqnarray}
where $A_\omega$ is the $\omega$-section of $A$. Moreover, the time average of the periodic measure defined by
\begin{eqnarray}\label{prop6b}
\bar \mu={1\over\tau}\int _0^\tau \mu _sds.
\end{eqnarray}
 is an invariant measure of $\Phi$ whose random factorisation is supported by $L^{\omega}$ defined in (\ref{prop6e}).
\end{thm}

\begin{proof}
It is obvious that $P$ is the marginal measure of $\mu_s$ on $(\Omega,{\mathcal F})$, so $\mu_s\in {\mathcal P}_P(\Omega\times \mathbb{\mathbb{X}})$.
To check (\ref{1.3}), first note for $t\in \mathbb{R}^{+}$,
$\bar\Theta(t)^{-1}(A)=\{(\omega, x): (\theta(t)\omega, \Phi(t,\omega)x)\in A\}.$
Then it is easy to see that for $t\in \mathbb{R}^{+}$
\begin{eqnarray}\label{zhao12}
(\bar\Theta_t^{-1}(A))_\omega
=
\Phi^{-1}(t,\omega)A_{\theta(t)\omega}.
\end{eqnarray}
Then (\ref{1.3}) follows a standard argument. 
Thus $\mu _s, s\in \mathbb{R}$ defined by (\ref{1.4}) is a periodic measure as claimed in the theorem.
To see $\bar \mu $ defined by (\ref{prop6b}) is an invariant measure, note for any
$A\in {\mathcal F}\otimes {\mathcal B}(\mathbb{\mathbb{X}})$ and $t\in \mathbb{R}^+$, by what we have proved for $\mu_s$,
\begin{eqnarray}\label{im1}
\bar \Theta (t)\bar \mu (A)={1\over\tau}\int _0^\tau \bar \Theta (t)\mu _s(A)ds= {1\over\tau}\int _0^\tau \mu _{t+s}(A)ds
=
{1\over\tau}\int _0^{\tau} \mu _{s}(A)ds
=
\bar \mu (A).
\end{eqnarray}
Thus $\bar \mu$ is an invariant measure. To see its support,
by (\ref{prop6b}), (\ref{1.4}) and Fubini's Theorem, for any $A\in \bar {\mathcal F}$,
\begin{eqnarray*}
 \bar \mu(A)
 ={1\over\tau}\int _0^\tau \mu_s(A)ds= {1\over\tau}\int _0^\tau \int_\Omega \delta_{Y(s,\theta(-s)\omega)}(A_{\omega})P(d\omega)ds=\int_\Omega {1\over\tau}\int _0^\tau  \delta_{Y(s,\theta(-s)\omega)}(A_{\omega})dsP(d\omega).
 \end{eqnarray*}
 This leads to its factorisation given by
 \begin{eqnarray*}
 (\bar\mu)_\omega={1\over\tau}\int _0^\tau  \delta_{Y(s,\theta(-s)\omega)}ds={1\over\tau}\int _0^\tau  \delta_{\phi(s,\omega)}ds,
 \end{eqnarray*}
which is supported by $L^{\omega}$.
\end{proof}
\begin{rmk}
For a random periodic path $Y$,
it is easy to see that the factorization of $\mu_s$ defined in Theorem \ref{zhao13} is
\begin{eqnarray}\label{10}
(\mu_s)_\omega=\delta_{Y(s,\theta(-s)\omega)},
\end{eqnarray}
and satisfies
\begin{eqnarray}
(\mu_{s+\tau})_\omega=(\mu_s)_\omega,\ \Phi(t,\omega)(\mu_s)_\omega=(\mu_{t+s})_{\theta(t)\omega}.
\end{eqnarray}
\end{rmk}

Now consider a Markovian cocycle random dynamical system $\Phi$ on a filtered dynamical system $(\Omega, {\mathcal F}, P, (\theta_t)_{t\in \mathbb{R}}, ({\mathcal F}_s^t)_{s\leq t})$, i.e. for any $s, t, u\in \mathbb{R}, s\leq t$, $\theta_u^{-1}{\mathcal F}_s^t={\mathcal F}_{s+u}^{t+u}$ and for any $t\in \mathbb{R}^+$, $\Phi(t,\cdot)$ is measurable with respect to ${\mathcal F}_0^t$. We also assume
the random periodic solution $Y(s)$ is adapted, that is to say that  for each $s\in \mathbb{R}$, $Y(s,\cdot)$ is measurable with respect to ${\mathcal F}_{-\infty}^s:=\vee_{r\leq s}{\mathcal F}_r^s$.

Denote the transition probability of Markovian process $\Phi(t,\omega)x$ on the Polish space $\mathbb{\mathbb{X}}$ with Borel $\sigma$-field ${\mathcal B}(\mathbb{\mathbb{X}})$ by (c.f.
Arnold \cite{ar}, Da Prato and Zabczyk \cite{da-prato}) 
$$P(t,x, B)=P(\{\omega:\Phi(t, \omega)x\in B\}), \ \ t\in \mathbb{R}^+,\  B\in {\mathcal B}(\mathbb{\mathbb{X}}),$$
and for any probability measure $\rho$ on  $(\mathbb{\mathbb{X}},{\mathcal B}(\mathbb{\mathbb{X}}))$, define
$$(P^*(t)\rho)(B)=\int_\mathbb{\mathbb{X}} P(t, x, B)\rho(dx),\ \ {\rm for \ any }\  t\in \mathbb{R}^+, \ B\in {\mathcal B}(\mathbb{\mathbb{X}}).$$
\begin{defi}\label{zhao14}
A measure function $\rho: \mathbb R\to {\mathcal P}(\mathbb{\mathbb{X}})$ is called a periodic measure of period $\tau$
on the phase space $(\mathbb{\mathbb{X}},{\mathcal B}(\mathbb{\mathbb{X}}))$ for the Markovian random dynamical systems $\Phi$ if it satisfies
\begin{eqnarray}\label{zhaoh61}
\rho_{s+\tau}=\rho_s\  \ {\rm and}\ \ \rho_{t+s}(B)=\int_\mathbb{\mathbb{X}} P(t, x, B)\rho_s(dx) \ \ s\in \mathbb{R},\ t\in \mathbb{R}^+.
\end{eqnarray}
It is called a periodic measure with minimal period $\tau$ if  $\tau>0$ if the smallest number such that (\ref{zhaoh61}) holds.
It is called an invariant measure if it also satisfies $\rho_s=\rho_0$ for all $s\in \mathbb{R}$, i.e. $\rho_0$ is an invariant measure for the Markovian random dynamical system $\Phi$ if
\begin{eqnarray}
\label{zhaoh622}
\rho_{0}=P^*(t)\rho_0, {\rm \ for \ all }\ t\in \mathbb{R}^+.
\end{eqnarray}
\end{defi}

\begin{rmk}
In \cite{Ha}, Has'minskii suggested that 
\begin{eqnarray}\label{zhao120}
P^*_{k\tau}\rho_s=\rho_{k\tau+s}=\rho_s.
\end{eqnarray}
It is easy to construct a counter example which satisfies (\ref{zhao120}), but not (\ref{zhaoh61}). For 
example we consider a RDS with two different periodic measures of the same period with non-overlapping supports. Then we can construct a new measure function $\{\rho_s\}_{s\in [0,\tau)}$ by choosing one periodic measure for certain time and another periodic measure for other time.
The measure function can be extended to all $s\in {\mathbb R}$ by 
imposing the periodicity in time.   Then the new measure function still satisfies (\ref{zhao120}), but does not satisfy (\ref{zhaoh61}). Certainly it does not make sense to say it is a periodic measure of the Markov semigroup as it is constructed from two different periodic measures. 
In fact, both conditions in the definition (\ref{zhaoh61}) are not redundant. 
  \end{rmk}
  
\begin{thm}\label{zhao31}
Assume the Markovian cocycle $\Phi: \mathbb{R}^+\times \Omega\times \mathbb{\mathbb{X}}\to \mathbb{\mathbb{X}}$ has an adapted random periodic path $Y: \mathbb{R}\times \Omega\to \mathbb{\mathbb{X}}$. Then
the measure function $\rho: \mathbb{R}\to {\mathcal{P}}(\mathbb{\mathbb{X}})$  defined by
\begin{eqnarray}\label{zhao11}
\rho_s:=E(\mu_s)_{\cdot}=E\delta_{Y(s,\theta(-s)\cdot)}=E\delta_{Y(s,\cdot)},
\end{eqnarray}
which is the law of the random periodic path $Y$,
is a periodic measure of $\Phi$ on
$(\mathbb{\mathbb{X}},{\mathcal{B}}(\mathbb{\mathbb{X}}))$.  Its time average $\bar \rho$ over a time interval of exactly one period
defined by
\begin{eqnarray}\label{zhao42}
\bar \rho={1\over\tau}\int _0^\tau \rho _sds,
\end{eqnarray}
is an invariant measure and satisfies that for any $B\in {\mathcal{B}}(\mathbb{\mathbb{X}})$, $t\in \mathbb{R}$
\begin{eqnarray}\label{zhao43}
&&
{\bar \rho}(B)\nonumber
\\
&=&E({1\over\tau} m\{s\in [0,\tau): Y(s,\cdot)\in B\})\nonumber
\\
&=&E({1\over\tau} m \{s\in [t,t+\tau): Y(s,\cdot)\in B\}).
\end{eqnarray}
\end{thm}
\begin{proof} Firstly it is easy to see from the definition of random periodic path  that for any $B\in {\mathcal{B}}(\mathbb{\mathbb{X}})$,
$\rho_{s+\tau}(B)
=\rho_s(B).$
Secondly, from (\ref{10}) we have
$(\mu_{t+s})_\omega=\delta_{Y(t+s,\theta(-t-s)\omega)}=\delta_{\Phi(t,\theta(-t)\omega)Y(s,\theta(-t-s)\omega)}.$
Therefore for any $B\in {\mathcal{B}}(\mathbb{\mathbb{X}})$, $t\in \mathbb{R}^+$, by measure preserving property of $\theta$, independency of $\Phi(t, \theta(s)\omega)$ and ${\mathcal{F}}_{-\infty}^s$,
\begin{eqnarray}\label{zhao41}
\rho_{t+s}(B)
&=&E\delta_{\Phi(t,\theta(s)\cdot)Y(s,\cdot)}(B)\nonumber
\\
&=&P(\{\omega:\Phi(t,\theta(s)\omega)Y(s,\omega)\in B\})\nonumber
\\
&=&
\int_\mathbb{\mathbb{X}} P(t,x,B)P(\omega: Y(s,\omega)\in dx)\nonumber\\
&=&\int_\mathbb{\mathbb{X}} P(t,x,B)\rho_s(dx)\nonumber\\
&=&P^*(t)\rho_s(B).
\end{eqnarray}
Therefore $\rho$ satisfies Definition \ref{zhao14} so is a periodic measure on $(\mathbb{\mathbb{X}},{\mathcal{B}}(\mathbb{\mathbb{X}}))$.
To prove the second part of the theorem, similar to the computation in (\ref{im1}), we have for any $t\in [0,\tau)$,
$\int_0^\tau\rho_{t+s}ds=\int_0^\tau\rho_{s}ds$,
and by using Fubini's Theorem,
\begin{eqnarray*}
\int_0^\tau P^*(t)\rho_sds= P^*(t)(\int_0^\tau\rho_{s}ds).
\end{eqnarray*}
It then follows easily that $\bar \rho$ is an invariant measure of $\Phi$ satisfying (\ref{zhaoh622}).
To prove the last part of the theorem, from (\ref{zhao42}), (\ref{zhao11}), and using Fubini's Theorem,  we know for any $B\in {\mathcal{B}}(\mathbb{\mathbb{X}})$,
\begin{eqnarray*}
{\bar \rho}(B)
&=& {1\over\tau}\int _0^\tau P(\omega: Y(s,\omega)\in B)ds\\
&=& {1\over\tau}\int _0^\tau (E \delta _{Y(s,\cdot)}(B))ds\\
&=& E({1\over\tau}\int _0^\tau  \delta _{Y(s,\cdot)}(B)ds)\\
&=& E({1\over\tau} m \{s\in [0,\tau): Y(s,\cdot)\in B\}).
\end{eqnarray*}
However, since ${\bar \rho}$ is an invariant measure, so from (\ref{zhao41}) we know that for any $t\in \mathbb{R}^+$
\begin{eqnarray}
\bar \rho(B)
&=&P^*(t)\bar\rho(B)\nonumber\\
&=&{1\over\tau}\int _0^\tau\rho_{t+s} (B)ds\nonumber\\
&=& E({1\over\tau} \int _t^{t+\tau}  \delta _{Y(s,\cdot)}(B)ds)\nonumber\\
&=& E({1\over\tau} m\{s\in [t,t+\tau): Y(s,\cdot)\in B\}).\label{3.17}
\end{eqnarray}
For $t\in \mathbb{R}^-$, it is easy to verify $P^*(-t)\rho_{s+t}=\rho_s$ and therefore
\begin{eqnarray*}
\bar \rho(B)
&=&P^*(-t)\bar\rho(B)\\
&=&{1\over\tau}\int _t^{t+\tau}P^*(-t)\rho_{s+t}  (B)ds\\
&=&{1\over\tau}\int _t^{t+\tau}\rho_{s}  (B)ds.
\end{eqnarray*}
Thus the result for $t\in \mathbb R^-$ follows a similar argument as (\ref{3.17}). In conclusion, we can see that (\ref{zhao43}) is true for any $t\in \mathbb{R}$.
\end{proof}
\vskip5pt

\begin{rmk}
The proof of $\bar \rho$ being an invariant measure does not depend
on $\rho_s$ being defined by (\ref{zhao11}). As long as $\rho_s,s\in {\mathbb R}$, is a periodic measure
of $\{P(t)\}_{t\geq 0}$, $\bar \rho$ defined by (\ref{zhao42}) is an invariant measure of $\{P(t)\}_{t\geq 0}$.
\end{rmk}

We observe that identity (\ref{zhao43}) says that the expected time spent inside a Borel set
 by the random periodic path over a time interval of exactly one period starting at any time is invariant, i.e. independent of the starting time.
This shows that the random periodicity of a random periodic path by means of invariant measures.
In the following we will establish the ergodic theory and the mean ergodic theory for periodic measures and 
random periodic paths.
They push (\ref{zhao43}) and the above observation much further in the case when a random periodic path exists.  
They say that on the long run, 
the average time that the random periodic path spends on a Borel set $B$ over one period is equal to $\bar \rho (B)$ both  
in law and in the long time average a.s.

\section{Poincar\'e sections and ergodicity with periodicity} 

We start to study the ergodicity of the random dynamical systems when periodicity exists. 
It is noted that the classical ergodic theorem dealing with invariant measures and stationary processes from Khas'minskii and Doob's
theorems fails to work in the stochastic periodic regime. Doob's classical method says that if the Markov transitional probability measures $P(t,x,\cdot)$, $x\in \mathbb{\mathbb{X}}$,
are mutually equivalent at a certain time $t_0>0$ ($t_0$-regular), then the invariant measure is strongly mixing 
and unique. Khas'minskii's theorem provides sufficient condition of verifying the regularity which says that if the Markovian semigroup is $t_0$-irreducible for certain $t_0>0$ and 
strong Feller at certain $t_1>0$, then the Markov semigroup is $(t_0+t_1)$-regular.

However, in a random periodic regime, if the periodic measure has a minimum period $\tau>0$, the 
invariant measure is not mixing and the $t_0$-regularity of the Markovian semigroup is no longer true any more.
This crucial assumption in Doob's Theorem excludes random periodic 
case automatically.  The irreducibility condition may not be always true on the whole space, in particular, if 
the support of $\rho_s$ is not the whole space for a given $s$, then for a nonempty open set $\Gamma$ lying outside of ${\rm supp}(\rho_s)$,
 $\Phi(s)x$ reaches $\Gamma$ with probability $0$ for any $x\in {\rm supp} (\rho_0)$.
 Even the irreducibility condition is satisfied, the strong Feller condition is a strict requirement which may not be satisfied in many situations 
 e.g. when the coppesponding second order differential operator, which is the infinitesimal generator of the Markovian semigroup in the 
 case of diffusions, is not strictly elliptic. 

Our idea here is to study, for any $s\in [0,\tau)$, 
the $\tau$-mesh discrete time random dynamical systems at integral multiples of the period $\Phi(k\tau,\omega): L_s\to L_s$, $k\in {\mathbb N}=\{0,1,2,\cdots\}$. Here 
$L_s={\rm supp}(\rho_s)$. For each $s\in [0,\tau)$, the measure $\rho_s$ on $L_s$ is an invariant measure with respect 
to $P(k\tau)$, $k\in \mathbb N$. This brings us back to the stationary regime. Then we are in the right set-up to discuss the irreducibility and mixing property 
of $\Phi(k\tau)$ on $L_s$. Then through the Markovian property, periodicity and the Chapman-Kolmogorov equation, we can obtain the ergodicity 
of the original random dynamical system if $\rho_s$ is a mixing invariant measure of the discrete Markov semigroup $P(k\tau)$. 

 We abstract the above idea to give the following definition first without assuming even the existence of periodic measures in the first place. 
 
 \begin{defi}\label{def2.22}
 The sets $L_s\subset \mathbb{\mathbb{X}},\ s\geq 0$, are called the Poincar$\acute{e}$ sections of the transition probability $P(t,\cdot,\cdot)$, $t\geq 0$, if $$L_{s+\tau}=L_s,
 $$ and
for any  $x\in L_{s}$, $t\geq 0$,
\begin{eqnarray}\label{zhao100c}
P(t,x, L_{s+t})=1.
\end{eqnarray}
\end{defi}

\begin{rmk}\label{zhao100m}
(i).
In fact, $L_s, s\geq 0$ can be extended to any $s\in \mathbb R$
 by the periodicity of $L_..$
 \\
 
 \noindent
(ii).
It is easy to see for each Poincar$\acute{e}$ section $L_s$, we have 
\begin{eqnarray*}
P(k\tau,x, L_{s})=1,\ \ {\rm for \ any \ } x\in L_{s}. 
\end{eqnarray*}
This means starting from $x\in L_s$, $
\Phi(k\tau,\omega)x$ returns to the set $L_s$ with probability 
one at any time being a multiple integral of the period.
This could be regarded as the Poincar\'e returning map property in the random regime, mirroring
the celebrated Poincar\'e mapping in the deterministic case. 
However, the map $\Phi(k\tau,\omega)$ does not have a fixed point on $L_s$. This is very different from 
the deterministic case.
\\

\noindent
(iii). 
It is worth pointing out that under the condition of existence of periodic measures, nontrivial Poincar\'e sections automatically 
exist. To see this, let $L_s:={\rm supp} (\rho_s)$. Then
for any $t\geq 0$
\begin{eqnarray}\label{20153a}
\int _{L_s}P(t,x,L_{t+s})\rho_s(dx)=\rho_{t+s}(L_{t+s})=1.
\end{eqnarray}
This, together with the fact that $0\leq P(t,x, L_{t+s})\leq 1$, implies that
\begin{eqnarray}\label{zhao100j}
P(t,x, L_{t+s})=1 {\  for \ }\rho_s-{almost\ all} \  x\in L_s, {\ for\  any }\ t\geq 0.
\end{eqnarray}

\noindent
(iv). Only from Definition \ref{def2.22}, the choice of the Poincar\'e sections may not be unique. For example in all cases, $L_s=\mathbb{\mathbb{X}}, s\geq 0$ is also a trivial choice of Poincar\'e sections satisfying Definition \ref{def2.22}. 
We will further add irreducibility condition below to guarantee a unique choice of the Poincar\'e sections up to a shift
(Lemma \ref{zhao311r}). But the irreducibility is not immediately needed in the following compactness theorem.

 \end{rmk}
 
 Recall a Markovian semigroup $\{P(t)\}_{t\geq 0}$, is said to be {\it stochastically continuous} (\cite{da-prato}) if 
$$
\lim_{t\to 0}P(t, x, B(x,\gamma))=1,\ {\rm for \ all} \ x\in \mathbb{\mathbb{X}}, \ \gamma>0.
$$
Denote by $B_b(\mathbb{\mathbb{X}})$, the space of all bounded Borel measurable functions on $\mathbb{\mathbb{X}}$, and $C_b(\mathbb{\mathbb{X}})$ the space of all bounded continuous
functions on $\mathbb{\mathbb{X}}$. For any $\phi \in B_b(\mathbb{\mathbb{X}})$, define
$$
P(t)\phi(x)=\int _\mathbb{X}P(t,x,dy)\phi(y), {\rm \ for }\ t\geq 0.
$$ 
Recall that the stochastically continuous 
semigroup $\{P(t)\}_{t\geq 0}$, is called a {\it Feller semigroup} if for any $\phi\in C_b(\mathbb{X})$,
we have $P(t)\phi\in C_b(\mathbb{X})$ for any $t\geq 0$. It is called a {\it strong Feller semigroup} at a time $t_0>0$ 
on a subset $\Gamma $ of $\mathbb{X}$ if for any $\phi\in B_b(\mathbb{X})$,
we have $P(t_0)\phi(x)|_{x\in \Gamma } \in C_b(\Gamma)$.

 Define now for any $\Gamma\in {\mathcal{B}}(\mathbb{X})$,
$$R_N(x, \Gamma):={1\over N}\sum_{k=1}^N P(k\tau,x, \Gamma),$$ 
and 
$$(R_N^*\nu)(\Gamma):=\int_\mathbb{X}R_N(x,\Gamma)\nu(dx),$$
for a measure $\nu\in {\mathcal P}(\mathbb{X})$.
Note if $\nu$ has a support in $L_0$, then
\begin{eqnarray*}
(R_N^*\nu)(L_0)={1\over N}\sum_{k=1}^N\int_\mathbb{X} P(k\tau,x, L_0)\nu(dx)=1.
\end{eqnarray*}
So ${\rm supp}(R_N^*\nu)\subset L_0$.

With the help of the Krylov-Bogoliubov theorem, we can prove the following existence theorem for a periodic measure.
  
\begin{thm}
Assume $L_s, s\in \mathbb R$ are Poincar\'e sections of Markovian semigroup $\{P(t)\}_{t\geq 0}$ and  $P(t)$ is a Feller semigroup on $L_0$. If for some $\nu\in {\mathcal P}(\mathbb{X})$ with its support in  $L_0$ and a subsequence $N_i$ with $N_i\to \infty$ as $i\to\infty$ such that
$$R_{N_i}^*\nu\to \rho_0,$$
weakly as $i\to \infty$. Define for any $\Gamma \in {\mathcal B}(\mathbb{X})$, if $s\geq 0$
\begin{eqnarray*}\rho_s(\Gamma)&=&\int_{L_0}P(s,x,\Gamma)\rho_0(dx),
\end{eqnarray*}
and if $ s<0$, 
\begin{eqnarray*}
\rho_s(\Gamma)&=& \rho_{s+k\tau}(\Gamma), 
\end{eqnarray*}
where $k$ is the smallest integer such that $s+k\tau \geq 0$.
Then $\rho _s, s\in {\mathbb R}$, is a periodic measure with respect to the semigroup $\{P(t)\}_{t\geq 0}$.  For each $s\in {\mathbb R}$, ${\rm supp}(\rho_s)\subset L_s$.
In particular $\rho_s(L_s)=1$.
\end{thm} 
\begin{proof} By the Krylov-Bogoliubov Theorem, it is easy to see that $\rho_0$ is an invariant measure of $P(k\tau), k=0,1,2,\cdots$, and $\rho_0(L_0)=1$. Thus ${\rm supp}(\rho_0)\subset L_0$ as $\rho_0$ is a probability measure.
From the definition of $\rho_s$, when $s\geq 0$, by (\ref{zhao100c}),
\begin{eqnarray*}
\rho_s(L_s)=\int_{L_0}P(s,x,L_s)\rho_0(dx)=\rho_0(L_0)=1.
\end{eqnarray*}
Similarly, ${\rm supp}(\rho_s)\subset L_s$. Thus when $s\geq 0$, by Chapman-Kolmogorov equation, Fubini's theorem and the fact that $\rho_0$ is the invariant measure of $P(\tau)$, for any $\Gamma\in \mathcal B(\mathbb{X})$,
\begin{eqnarray*}
\rho_{s+\tau}(\Gamma)
&=& \int_{L_0}\int_\mathbb{X} P(s,y,\Gamma)P(\tau,x, dy)\rho_0(dx)\\
&=&\int_\mathbb{X} P(s,y,\Gamma)\int_{L_0}P(\tau,x, dy)\rho_0(dx)\\
&=&\int_\mathbb{X} P(s,y,\Gamma)\rho_0(dy)\\
&=&\rho_s(\Gamma).
\end{eqnarray*}
Moreover, for any $t\geq 0$, $s\geq 0$, by a similar argument as above,
\begin{eqnarray*}
(P^*(t)\rho_s)(\Gamma)
&=&\int_{\mathbb{X}}P(t,x,\Gamma)\int_{L_0}P(s,y,dx)\rho_0(dy)\\
&=&\int_{L_0}\int_\mathbb{X} P(t,x,\Gamma)P(s,y,dx)\rho_0(dy)\\
&=&\int_{L_0} P(t+s,y,\Gamma)\rho_0(dy)\\
&=&\rho_{t+s}(\Gamma).
\end{eqnarray*}
That is to say $\rho_s, s\geq 0$ is the periodic measure of the transition semigroup $\{P(t)\}_{t\geq 0}$.
For $s<0$, it is obvious to verify the result.
\end{proof}

This theorem could be regarded as the extension of Krylov-Bogoliubov theorem to the periodic measure case. Though the theorem looks very different from the Poincar$\rm \acute{ e}$-Bendixson theorem in the first sight, but in spirit it is indeed like the Poincar\'e-Bendixson theorem as a random counterpart in the level of measures.
Though the Poincar\'e map does not have a fixed point in the pathwise sense, but 
$$
P^*(k\tau): {\mathcal P}(L_s)\to {\mathcal P}(L_s)
$$
 has a fixed point $\rho_s\in {\mathcal P}(L_s)$
for all $s\in {\mathbb R}$. All these invariant measures of $P(k\tau)$ together form a periodic measure.

Now we start to consider ergodicity. Recall $\bar\rho={1\over \tau}\int_0^\tau \rho_s ds$ as the invariant measure. 
Consider a set of finite sequence of real numbers ${\mathbb I}=\{t_1,t_2,\cdots,t_n\}$, $t_1<t_2<\cdots <t_n$ and by the Chapman-Kolmogorov equation and a standard procedure (\cite{da-prato}), we can construct $\{P_{\mathbb I}^{\bar \rho}, {\mathbb I}$ as a set of sequences of distinct real numbers$\}$ is a consistent family of finite dimensional
 distributions, where 
 $$P_{\mathbb I}^{\bar \rho}(\mathbb X\times\cdots \mathbb X\times A\times\cdots \times\mathbb X)={\bar \rho} (A).$$
By the Kolmogorov extension theorem, there exists a unique probability measure $P^{\bar \rho}$ on 
$(\Omega ^*, {\mathcal F}^*)=(\mathbb{X}^{\mathbb R}, {\mathcal B}(\mathbb{X}^{\mathbb R}))$ with a family of finite-dimensional distributions
$\{P_{\mathbb I}^{\bar \rho}\}_{\mathbb I}$. For any $\omega^*\in \Omega ^*$, denote its canonical process 
by 
$
W_t(\omega^*)=\omega^*(t),
$
 which is a Markovian process and a
measurably invertible map $\theta^*:{\mathbb R}\times \Omega^*\to \Omega^*$ by
$(\theta _t^*\omega^*)(s)=\omega^*(t+s), \ t,s\in {\mathbb R}.$
It follows that $(\Omega ^*, {\mathcal F}^*, \theta _t^*, P^{\bar \rho})$ defines a dynamical system, which is called 
the canonical dynamical system
associated with the semigroup $P_t, t\geq 0$ and invariant measure $\bar\rho$. It is well known that if $P_t, t\geq 0$ is
stochastically continuous, 
then the linear transformation operator $U_t:{\mathcal H}_{\mathbb C}^{\bar\rho}\to {\mathcal H}_{\mathbb C}^{\bar\rho}$, where
${\mathcal H}_{\mathbb C}^{\bar\rho}=L^2_{\mathbb C}(\Omega^*,{\mathcal F}^*,P^{\bar\rho})$ defined by 
\begin{eqnarray}
\label{zhao20157a}
U_t\xi(\omega^*)=\xi(\theta_t^*\omega^*), \ \xi \in {\mathcal H}^{\bar \rho}_{\mathbb C}, \ \omega^*\in \Omega^*, t\in {\mathbb R},
\end{eqnarray}
is continuous in $t$, and $(\Omega ^*, {\mathcal F}^*, \theta _t^*, P^{\bar \rho})$ is a continuous metric dynamical system. 
The invariant measure $\bar\rho$ 
is called ergodic if $(\Omega ^*, {\mathcal F}^*, \theta _t^*, P^{\bar \rho})$ is ergodic i.e.
\begin{eqnarray*}
\lim_{T\to\infty}{1\over T}\int_0^T P^{\bar\rho}(\theta_{-t}^*A\cap B)dt=P^{\bar\rho}(A)P^{\bar\rho}(B), \ \ {\rm for \ any}\ A,B\in {\mathcal F}^*.
\end{eqnarray*}
We say that the periodic measure $\{\rho_t\}_{t\in {\mathbb R}}$, is ergodic if its average $\bar\rho$ as an invariant measure is ergodic. 
Also recall that an invariant measure $\rho$ 
is called weakly mixing if $(\Omega ^*, {\mathcal F}^*, \theta _t^*, P^{\rho})$ is weakly mixing i.e.
there is a set $I\subset [0,\infty)$ of relative measure 1 such that
\begin{eqnarray*}
\lim_{t\to\infty, t\in I} P^{\bar\rho}(\theta_{-t}^*A\cap B)=P^{\bar\rho}(A)P^{\bar\rho}(B), \ \ {\rm for \ any}\ A,B\in {\mathcal F}^*.
\end{eqnarray*}
The ergodicity and mixing property of discrete random dynamical systems, which will also be needed in this paper, can also be defined similarly by replacing the integral by summation and 
$\lim\limits _{t\to\infty, t\in I}$ by the limit along the discrete time sequence respectively. 

It is well-known that the following statements are equivalent (c.f. Theorem 3.2.4, \cite{da-prato}):

(i) $\bar\rho$ is ergodic;

(ii) if $U_t\xi=\xi$ for all $t\in {\mathbb R}^+$, then $\xi$ is constant;

(iii) if $P(t)\phi=\phi$ for all $t\in {\mathbb R}^+$, then $\phi$ is a constant; 

(iv) if 
a set $\Gamma\in {\mathcal B}(\mathbb{X})$ satisfies for all $t\in {\mathbb R}^+$
\begin{eqnarray*}
P_tI_{\Gamma}=I_{\Gamma},\ \bar\rho-a.e.
\end{eqnarray*}
then either $\bar\rho(\Gamma)=0$ or $\bar\rho(\Gamma)=1$;

(v) for any $\Gamma\in {\mathcal B}(\mathbb{X})$, $\lim\limits_{T\to\infty}{1\over T}\int_0^TP(s,x,\Gamma)ds\to \bar\rho (\Gamma), \ {\rm in }\ L^2({\mathbb X},\bar\rho(dx)).
$

\noindent
Moreover, the following statements are also equivalent:

(vi) $\rho$ is weakly mixing;

(vii) if $U_t\xi={\rm e}^{i\lambda t} \xi$ for all $t\in {\mathbb R}^+$, $\lambda $ is a real number,  then $\lambda=0$ and $\xi$ is constant;

(viii)
if $P(t)\phi={\rm e}^{i\lambda t} \phi $ for all $t\in {\mathbb R}^+$, $\lambda $ is a real number,  then $\lambda=0$ and $\phi$ is a constant;

(ix)
there exists $I\subset [0,\infty)$ of relative measure 1 such that 
$$
\lim\limits_{t\to \infty, t\in I} P(t,x,-)\to \rho .
$$

The equivalence of (vi) and (vii) is the Koopman-von Neumann Theorem.
From equivalence of (vi) and (ix) in the above, it is easy to see there is no way one can establish the mixing property for $\bar \rho$ in the regime of 
random periodicity unless it is degenerated to the stationary case.

Now we assume a periodic measure $\{\rho_s\}_{s\in {\mathbb R}}$ exists with $L_s={\rm supp}(\rho_s)$. Set 
\begin{eqnarray}\label{fengc2}
L:=\bigcup\{L_s: s\in [0,\tau)\}.
\end{eqnarray} 
Then it is easy to see that 
\begin{eqnarray}
\bar \rho (L)={1\over \tau}\int _0^{\tau} \rho_s(L)ds=1.
\end{eqnarray}
This implies that ${\rm supp}(\bar\rho)\subset L$. 
Moreover, note $\bar\rho(B)=1$ iff $\rho_s(B)=1$ for almost all $s\in [0,\tau)$. So $L_s\subset {\rm supp}(\bar\rho)$. Thus $
L=\bigcup\limits _{0\leq s\leq \tau} L_s= {\rm supp}(\bar\rho)$.

We first prove a simple but useful lemma. For this we consider
\\

\noindent
{\bf Condition A}: {\it The Markovian cocycle $\Phi: \mathbb{R}^+\times \Omega\times \mathbb{X}\to \mathbb{X}$  has a periodic measure $\rho: \mathbb{R}\to {\mathcal{P}
}(\mathbb{X})$ and for any $\Gamma \in {\mathcal{B}}(\mathbb{X})$, we have when $N\to \infty$,
\begin{eqnarray}
\int_\mathbb{X}|\int_{0}^{\tau} ({1\over N}\sum\limits _{k=0}^{N-1}P(s+k\tau, y, \Gamma)-\rho_s(\Gamma))ds|\bar \rho(dy)\to 0,
\end{eqnarray}
where $\bar \rho={1\over \tau}\int _0^{\tau}\rho_sds$.
}
\\

\begin{lem}\label{fengc4}
Assume the Markovian semigroup $P(t)$ is stochastically continuous. Then 
 the invariant measure $\bar\rho$ is ergodic if and only if Condition A holds. Moreover, in this case
$L$ defined by (\ref{fengc2}) is the unique set (up to a $\bar \rho$-measure 0 set) with positive $\bar \rho$-measure satisfying $P(t,x,L)=I_L(x)$.
\end{lem}
\begin{proof} First assume Condition A holds. For any $\Gamma \in {\mathcal B}(\mathbb{X})$, if 
$$(P(t)I_{\Gamma})(\cdot)=P(t,\cdot,\Gamma )=I_{\Gamma}(\cdot),\ \bar\rho-a.e.,$$
then it turns out from Condition A that  
\begin{eqnarray*}
\int_\mathbb{X}|{1\over \tau}\int_{k\tau}^{(k+1)\tau} I_{\Gamma }(y)ds-\bar\rho(\Gamma )|\bar \rho(dy)=\int_\mathbb{X}| 
I_{\Gamma }(y)-\bar\rho(\Gamma )|\bar \rho(dy)= 0,
\end{eqnarray*}
so $$\bar\rho(\Gamma )=I_{\Gamma }(y),\ \bar\rho-a.e.$$
This implies that $I_{\Gamma}(y)$ is a constant for $\bar\rho$-a.e. $y\in \mathbb{X}$. Thus $$\bar\rho(\Gamma )=0\ \ {\rm or} \ \  1.$$
By Theorem 3.2.4 in \cite{da-prato}, $\bar\rho$ is ergodic. Moreover, from the fact that ${\rm supp}(\bar\rho)\subset L$, it is easy to see that in the case 
$\bar\rho(\Gamma)=1$, $\Gamma=L$ up to a $\bar\rho$-measure $0$ set.  The last claim is proved.

Conversely, assume $\bar\rho$ is ergodic. Then
${1\over T}\int_0^T P(s,x,\Gamma)ds\to \bar\rho(\Gamma)$ in $L^2({\mathbb X}, \bar\rho(dx))$. 
Thus ${1\over N\tau}\sum \limits _{k=0}^{N-1}\int_0^{\tau} P(s+k\tau,x,\Gamma)ds\to \bar\rho(\Gamma)$ in $L^2({\mathbb X},\bar\rho(dx))$.
Then Condition A follows from above and Cauchy-Schwarz inequality. 
 \end{proof}
 
With this lemma, we only need to verify Condition A in order to prove the ergodicity for an invariant measure generated 
by periodic measures. 

\begin{defi} The $\tau$-periodic measure  $\{\rho_s\}_{s\in {\mathbb R}}$ is called to be PS-ergodic (PS-mixing) if 
for each $s\in [0,\tau)$, $\rho_s$ as the invariant measure of 
 the $\tau$-mesh discrete Markovian semigroup $\{P(k\tau)\}_{k\in {\mathbb N}}$, at integral multiples of the period on 
the  Poincar\'e section 
 $L_s$, is ergodic (mixing). 
 \end{defi}

\begin{thm}\label{zhao311s}
Let the Markovian semigroup $P(t)$ be stochastically continuous and have a $\tau$-periodic measure $\{\rho_s\}_{s\in {\mathbb R}}$. 
 Assume  $\{\rho_s\}_{s\in {\mathbb R}}$ is PS-ergodic.  Then 
Condition A is satisfied and the invariant measure $\bar\rho$ is ergodic. In particular, $\Gamma=L$ is the unique set with positive $\bar \rho$-measure 
 satisfying 
$P(t,x,\Gamma)=I_{\Gamma}(x)$  for any $t\ge 0$.
Moreover, if $\tau>0$ is the minimum period of the periodic measure, then $L_{s_1}\cap L_{s_2}=\emptyset $ when 
$s_1,s_2\in [0,\tau), s_1\ne s_2$. 
\end{thm}
\begin{proof}
According to Theorem 3.4.1 in \cite{da-prato}, as for any fixed $s\in [0,\tau)$,  $\rho_s$ as the invariant measure of $P(k\tau)|_{L_s}, k\in {\mathbb N}$, is ergodic, so for any $\phi\in L^2(L_s,\rho_s)$, 
we have as $N\to \infty$,
\begin{eqnarray*}
{1\over N}\sum \limits _{k=0}^{N-1}P(k\tau)\phi(\cdot)\to <\phi ,1>_{L^2(L_s,\rho_s)},\ \ {\rm \ in } \ L^2(L_s,\rho_s).
\end{eqnarray*}
Now consider $\phi (\cdot)=P(t,\cdot,\Gamma)$ for an arbitrarily given $\Gamma \in {\mathcal B}(\mathbb{X})$. Note 
$P(k\tau)\phi(\cdot)=P(t+k\tau,\cdot,\Gamma)$ and $<\phi ,1>_{L^2(L_s,\rho_s)}=\rho_{t+s}(\Gamma).$ Thus as $N\to \infty$
\begin{eqnarray}\label{zhao311t}
{1\over N}\sum \limits _{k=0}^{N-1}P(t+k\tau,\cdot,\Gamma)\to \rho_{t+s}(\Gamma),\ \ {\rm \ in } \ L^2(L_s,\rho_s).
\end{eqnarray}
It then follows by applying Fubini theorem, Jensen's inequality and Lebesgue's dominated 
convergence theorem that
\begin{eqnarray*}
&&\int_\mathbb{X}|\int_{0}^{\tau} ({1\over N} \sum \limits _{k=1}^{N-1}P(t+k\tau, x, \Gamma)-\rho_t(\Gamma))dt|\bar \rho(dx)\\
&=&{1\over \tau}\int_0^\tau \int_\mathbb{X}|\int_{0}^{\tau}[{1\over N}\sum \limits _{k=1}^{N-1}P(t+k\tau, x, \Gamma)-\rho_{t+s}(\Gamma)]dt|\rho_{s}(dx)ds\\
&\leq&{1\over \tau}\int_0^\tau \int_{0}^{\tau} \int_\mathbb{X}|{1\over N}\sum \limits _{k=1}^{N-1}P(t+k\tau, x, \Gamma)-\rho_{t+s}(\Gamma)|\rho_{s}(dx)dtds\\
&\leq &{1\over \tau}\int_0^\tau \int_{0}^{\tau} \left [\int_\mathbb{X}|{1\over N}\sum \limits _{k=1}^{N-1}
P(t+k\tau, x, \Gamma)-\rho_{t+s}(\Gamma)|^2\rho_{s}(dx)\right ]^{1\over 2}
dtds\\
&\to& 0,
\end{eqnarray*}
as $k\to \infty$. Thus Condition A holds and the other results of the first part of the theorem follow. 

To prove the last result, from the PS-ergodicity of the periodic measure, we know that $
{1\over N}\sum \limits_{k=0}^{N-1}P(k\tau,x,\Gamma) \to \rho_{s}(\Gamma)$ in 
$L^2(\mathbb X, \rho_s(dx))$, for any $\Gamma\in {\mathcal B}(\mathbb X)$. 
So there exists a subsequence such that along the subsequence, the above convergence 
holds for $\rho_s$-a.e. $x$.
As $\tau>0$ is a minimum period, so for any $s_1,s_2\in [0,\tau), s_1\ne s_2$, 
$\rho_{s_1}\ne \rho_{s_2}$. Let $\Gamma \in {\mathcal B}(\mathbb X)$ be such that $\rho_{s_1}(\Gamma)\ne \rho_{s_2}(\Gamma)$. For $\rho_{s_1}, \rho_{s_2}$, there exists a common
subsequence $N_m\to \infty$ as $m\to \infty$ such that for any $\Gamma \in {\mathcal B}({\mathbb X})$, ${1\over N_m}\sum \limits_{k=0}^{N_m-1}P(k\tau,x,\Gamma) \to \rho_{s_1}(\Gamma)$ for $\rho_{s_1}$-a.e.
$x$ and ${1\over N_m}\sum \limits_{k=0}^{N_m-1}P(k\tau,x,\Gamma) \to \rho_{s_2}(\Gamma)$ for $\rho_{s_2}$-a.e.
$x$.
Set
\begin{eqnarray*}
A&=&\{x\in \mathbb{X}: {1\over N_m}\sum \limits_{k=1}^{N_m-1}P(k\tau,x,\Gamma) \to \rho_{s_1}(\Gamma)\},\\
B&=&\{x\in \mathbb{X}: {1\over N_m}\sum \limits_{k=1}^{N_m-1}P(k\tau,x,\Gamma) \to \rho_{s_2}(\Gamma)\}.\\
\end{eqnarray*}
So $\rho_{s_1}(A)=1$ and $\rho_{s_2}(B)=1$.  
But it is clear that $A\cap B= \emptyset$. Thus the last claim of the theorem is asserted.
\end{proof} 

Now we study the irreducibility condition. For $0\leq s<\tau$, consider

\begin{defi}
{\bf  (The $k_s\tau$-irreducibility condition on a Poincar$\rm \acute{\bf e}$ section $L_{s}$)}:  For a fixed $s\in [0,\tau)$,
 if there exists $k_s\in {\mathbb N}\setminus \{0\}$, such that
for an arbitrary 
nonempty relatively open set $\Gamma\subset L_{s}$, we have
\begin{eqnarray}
P(k_s\tau,x,\Gamma)>0, {\rm \ for \ } \rho_s-{\rm a.e. \ } x\in L_s,
\end{eqnarray}
then we call the Markovian semigroup $\{P(t)\}_{t\geq 0}$, is $k_s\tau$-irreducible on the Poincar\'e section $L_{s}$.
If for a certain map $s\mapsto k_s\in {\mathbb N}\setminus \{0\}, s\in [0,\tau)$, the semigroup is $k_s\tau$-irreducible for each $s\in [0,\tau)$, 
then we call the Markovian semigroup is $k_s\tau, s\in [0,\tau)$, irreducible on Poincar\'e sections $L_{s}, s\in [0,\tau)$. 
\end{defi}

\begin{defi}
{\bf  (The $k_s\tau$-regularity on a Poincar$\rm \acute{\bf e}$ section $L_{s}$)}:  A Markovian semigroup $\{P(t)\}_{t\geq 0}$, is said to be $t_0$-regular 
if all transitional probability measures $P(t_0,x,\cdot), x\in \mathbb{X}$, are mutually equivalent. For a fixed $s\in [0,\tau)$, it is said to be $k_s\tau$-regular for 
a certain $k_s\in {\mathbb N}\setminus \{0\}$ on a Poincar\'e section $L_s$, if all transitional probability measures $P(k_s\tau,x,\cdot), x\in L_s$, are mutually equivalent.
If for a certain map $s\mapsto k_s\in {\mathbb N}\setminus \{0\}, s\in [0,\tau)$, the semigroup is $k_s\tau$-regular on $L_s$ for each $s\in [0,\tau)$, 
then we call the Markovian semigroup is $k_s\tau, s\in [0,\tau)$, regular on Poincar\'e sections $L_{s}, s\in [0,\tau)$. 
\end{defi}

\begin{lem}\label{zhao311r}
Assume the Markovian semigroup $\{P(t)\}_{t\geq 0}$ has Poincar\'e sections $\{L_s\}_{s\in {\mathbb R}}$and periodic measure$\{\rho_s\}_{s\in {\mathbb R}}$ with ${\rm supp}(\rho_s)\subset L_s$. If $P(t)$ satisfies the $k_{s_0}\tau$-irreducibility condition on $L_{s_0}$ for some 
$k_{s_0}\in \mathbb N\setminus \{0\}$, then $L_{s_0}={\rm supp}(\rho_{s_0})$. Moreover, 
if  the semigroup satisfies the $k_s\tau$-irreducibility condition on the Poincar\'e sections $L_{s}$ for all $s\in [0,\tau)$, then $L_{s}={\rm supp}(\rho_{s})$ for any $s\in {\mathbb R}$.
\end{lem}
\begin{proof} By the $k_{s_0}\tau$-irreducibility condition on a Poincar\'e section $L_{s_0}$, 
we know that 
there exists $k_{s_0}\in {\mathbb N}\setminus \{0\}$ such that
for an arbitrary 
nonempty relatively open set $\Gamma\subset L_{{s_0}}$, we have
\begin{eqnarray*}
P(k_0\tau,x,\Gamma)>0, {\rm \ for \ } \rho_{s_0}-{\rm a.e. \ } x\in L_{s_0}.
\end{eqnarray*}
So for this $\Gamma$, we have 
\begin{eqnarray*}
\rho_{s_0}(\Gamma)=\int_{L_{s_0}} P(k_0\tau,x,\Gamma)\rho_{s_0}(dx)>0.
\end{eqnarray*}
Thus $L_{s_0}={\rm supp}(\rho_{s_0})$. The last claim follows easily from the above.
\end{proof}

\begin{rmk}
Under the irreducible conditions on Poincar\'e sections, it is easy to know that for any fixed $s\in [0,\tau)$ and  any open set $\Gamma _s\subset L_s={\rm supp}(\rho_s)$ with $\rho_s(L_s\setminus {\bar \Gamma }_s)>0$, we have
for any $x\in L_s$,
\begin{eqnarray*}
P(k\tau,x,\bar \Gamma _s)<1.
\end{eqnarray*}
This suggests that $\Gamma_s$ does not satisfy the requirement being a Poincar\'e section at time $s$. Thus, $L_s={\rm supp}(\rho_s)$, $s\in [0,\tau)$
are minimal Poincar\'e sections. 
\end{rmk}

\begin{thm}\label{thm2.23}
Let the Markovian semigroup $P(t)$ be stochastically continuous and have a $\tau$-periodic measure $\{\rho_s\}_{s\in {\mathbb R}}$. Denote $L_s={\rm supp}(\rho_s)$
 and $L=\bigcup_{0\leq s<\tau}L_s$. Assume the semigroup 
  is $k_s\tau$-regular, $s\in [0,\tau)$, on Poincar\'e sections for certain map $s\mapsto k_s\in {\mathbb N}\setminus \{0\}$, $s\in [0,\tau)$.  
 Then the periodic measure is PS-mixing and thus ergodic. 
\end{thm}
\begin{proof} 
Note first that $\rho_s$ is an invariant measure w.r.t. $P(k\tau),$ for any $k\in \mathbb N$. 
Due to the $k_s\tau$-regularity assumption, Doob's theorem (\cite{doob}) can be then applied to the discrete semigroup on the Poincar\'e section $L_s$ so the invariant measure $\rho_s$ of 
$\{P(k\tau)\}_{k\in {\mathbb N}}$, is ergodic on $L_{s}$ and for any $x\in L_{s}$, $\Gamma\in {\mathcal{B}}(\mathbb{X}),$
 \begin{eqnarray}\label{zhao100u}
P(k\tau,x,\Gamma)\to \rho_{s}(\Gamma),\ \ 
{\rm as \ } k\to \infty.
\end{eqnarray}
To see this, first note that Doob's theorem implies that (\ref{zhao100u}) holds for any $\Gamma \in {\mathcal B}(\mathbb{X})\cap L_s$. But (\ref{zhao100u}) is true for any
$\Gamma\in {\mathcal{B}}(\mathbb{X})$, as for any $x\in L_s$, $P(k\tau,x,\Gamma)=P(k\tau,x,\Gamma \cap L_s)$ and
$\rho_s(\Gamma)=\rho_s(\Gamma \cap L_s)$ since ${\rm supp}(\rho_s)=L_s$.
Therefore $P(k\tau,x,\cdot)\to \rho_{s}(\cdot)$ weakly by Proposition 2.4 in \cite{watanabe}. This implies that 
 the periodic measure is PS-mixing. Thus it is PS-ergodic and thus ergodic.
\end{proof}
The regularity of the semigroup condition can be checked.
\begin{lem}
 Assume the Markovian semigroup $P(t), t\geq 0,$ is $k_s\tau$-irreducible, $s\in [0,\tau)$, on  the Poincar\'e sections 
 $L_s={\rm supp}(\rho_s), s\in [0,\tau)$
 for certain map $s\mapsto k_s\in {\mathbb N}\setminus 
 \{0\}$, and strong Feller at $k_s^*\tau$ on $L_s$ for each $s\in [0,\tau)$, where $s\mapsto k_s^*\in {\mathbb N}\setminus 
 \{0\}$ is a certain map.  
  Then the semigroup is $(k_s+k_s^*)\tau$-regular on the Poincar\'e sections.
\end{lem}
\begin{proof} The proof is done by a similar proof as the one of Khas'minskii's theorem (\cite{Ha}) on each Poincar\'e section.
\end{proof}

\section
{Random periodic verses stationary:  sufficient-necessary conditions}

It is not a trivial task to check whether or not the minimum period of a random periodic solution is strictly positive. 
In this section, we will develop some equivalent sufficient and necessary 
conditions in four different notions. In particular, we will characterise it with an analytic assumption
that the infinitesimal generator of the corresponding Markov semigroup of the random dynamical 
system
 has infinitely many simple eigenvalues $\{{2m\pi\over \tau}i\}_{m\in {\mathbb Z}}$,
and no other eigenvalues on the imaginary axis. 

First note it is evident that if the invariant measure $\bar\rho$ is ergodic, 
and there exists a set $I\subset [0,\tau)$ with positive Lebesgue measure such that for each $s\in I$, $\rho_s$ is not ergodic with respect to $P(t), t\geq 0$,
then $\bar\rho\ne \rho_s$ for any $s\in I$. In this case, the periodic measure is not degenerated to an invariant measure.
 In the following we will mainly consider the case when the periodic measure $\{\rho_s\}_{s\in {\mathbb R}}$ is PS-ergodic.
 
Recall first the following standard definition.

\begin{defi}
Let $U_t : {\mathcal H}_{\mathbb C}^{\bar\rho}\to {\mathcal H}_{\mathbb C}^{\bar\rho}$ be the 
transformation operator defined by (\ref{zhao20157a}). A measurable function $\alpha: \Omega ^*\to [0,2\pi)$ is said to be an angle variable for the canonical 
dynamical system $(\Omega^*, {\mathcal F}^*, (\theta ^*(t))_{t\in {\mathbb R}}, 
P^{\bar \rho})$, if there exists a 
constant $\lambda \in {\mathbb R}$ such that for every $t\in {\mathbb R}$, 
\begin{eqnarray}\label{zhao20157b}
U_t\alpha =\lambda t+\alpha\  ({\rm mod} \ 2\pi), P^{\bar\rho}-a.s.
\end{eqnarray} 
\end{defi}

\begin{rmk}\label{zhao20157d}
(i) If $\alpha$ is an angle variable with $\lambda$, then $U_t e^{i\alpha}=e^{i\lambda t}e^{i\alpha}$, so $e^{i\lambda t}$ is an eigenvalue of $U_t$ and $\xi=e^{i\alpha}$ is the corresponding eigenvector in ${\mathcal H}_{\mathbb C}^{\bar\rho}$.

(ii) The following results in this paragraph are also well-known. We summarise them here as they are needed.   
As $U_t$ is a unitary operator and $U_t^*=U_{-t}$, so according to Stone's theorem, the infinitesimal operator of $U_t, t\in {\mathbb R}$, is of the form 
$i{\mathcal A}$, where ${\mathcal A}$ is a self-adjoint operator acting on ${\mathcal H}_{\mathbb  C}^{\bar\rho}$. The operator ${\mathcal A}$ is called the infinitesimal 
generator of the canonical dynamical system 
$(\Omega^*, {\mathcal F}^*, (\theta ^*(t))_{t\in {\mathbb R}}, P^{\bar \rho})$. Assume there exist $\lambda 
\in {\mathbb R}$ and $\xi \in {\mathcal H}_{\mathbb C}^{\bar\rho}\cap D({\mathcal A})$ such that 
\begin{eqnarray}\label{sept1}
{\mathcal A}\xi=\lambda \xi.
\end{eqnarray}
Then
\begin{eqnarray}\label{zhao20157c}
U_t\xi={\rm e}^{i{\mathcal A}t}\xi= {\rm e}^{i\lambda t}\xi.
\end{eqnarray}
It then follows that $U_t|\xi|=|\xi|$, $t\in {\mathbb R}$. So if $\bar \rho$ is ergodic, 
then $\xi$ is a constant and we can assume that $|\xi|=1.$ Consequently $\xi={\rm e}^{i\alpha}$, where 
$\alpha $ is a real valued random variable with values on $[0,2\pi)$. From (\ref{zhao20157c}), we know
that $\alpha$ is an angle variable satisfying (\ref{zhao20157b}). 
\end{rmk}

Recall the Koopman-von Neumann theorem which says that $\bar \rho$ is weakly mixing if and only if any angle variable is constant and the operator 
${\mathcal A}$ has only one eigenvalue $0$. Moreover, $0$ is a simple eigenvalue of ${\mathcal A}$. 

Note that the semigroup $P(t)$ is a map from $L^2(\mathbb{X},d\bar\rho)$ to $L^2(\mathbb{X},d\bar\rho)$. 
Recall the following well-known result (Theorem 3.2.1 in \cite{da-prato}): there exist, 
$\xi\in {\mathcal H}_{\mathbb C}^{\bar\rho}$, $\gamma \in {\mathbb C}$ with $|\gamma |=1$ such that 
$$
U_t\xi=\gamma \xi,
$$
 iff there exist $\phi\in L^2(\mathbb{X},d\bar\rho)$, $\gamma \in {\mathbb C}$ with $|\gamma |=1$  such that 
$$
P(t)\phi=\gamma \phi
$$
 and $\xi(\omega^*)=\phi(\omega^*(0))$. That is to say that all the eigenvalues of $P(t)$ on the unit circle agree 
with all the eigenvalues of the $U_t$. This will help in the proof of the next theorem to identify the spectra of  semigroup $U_t$ 
on the space of square integrable functions on the path space  $\mathbb{X}^{\mathbb R}$ to the spectra on the unit circle 
of semigroup $P(t)$ on the space of square integrable functions on the phase space $\mathbb{X}$. 

It is worth noting that the spectral analysis of the latter is easier to handle than the former one.
Moreover, the spectral structure of the latter is richer than that of the former one. This extra information of
the spectra of the semigroup gives more information about the dynamics of the Markov random 
dynamical system, e.g. 
mixing property and convergence rate of the transitional probability to the invariant measure in the stationary case.
We will prove in the next subsection that spectral gap of the 
semigroup on $L^2(L_s,\rho_s)$ for all $s\in {\mathbb R}$ 
leads to the PS-mixingness of $\{\rho_s\}_{s\in {\mathbb R}}$ 
and the mixing rate is given by the spectral gap. 
 
Moreover, the spectra of the semigroup $P(t)$ can be 
analysed by studying the spectra of its infinitesimal generator. 
 Recall the definition of the infinitesimal generator ${\mathcal L}$ 
 of the semigroup $P(t): L^2(\mathbb{X},d\bar\rho)\to L^2(\mathbb{X},d\bar\rho)$ given  by
\begin{eqnarray}\label{zhao300a}
{\mathcal L}\phi =\lim_{t\to 0+} {P(t)\phi-\phi\over t},
\end{eqnarray}
for all $\phi\in D({\mathcal L})$, where
$$
D({\mathcal L}):=\{\phi\in L^2(\mathbb{X},d\bar\rho): \lim_{t\to 0+} {P(t)\phi-\phi\over t} {\rm \ exists\ in\ }  L^2(\mathbb{X},d\bar\rho)\}.
$$

Following Theorem \ref{zhao311s}, we are now ready to 
prove the following theorem. 

\begin{thm}\label{zhao20157f}
Assume the transition probability is stochastically continuous and 
has a periodic measure $\{\rho_s\}_{s\in {\mathbb R}}$ of period $\tau$. Assume 
the $\tau$-periodic measure is PS-mixing. 
Then one of the following three cases happens:

\item{\underline {Case (i)}}.  The period $\tau$ is the smallest number such that (\ref{zhaoh61}) holds.

\item{\underline {Case (ii)}}. There exist $k\in \mathbb N\setminus \{0,1\}$, $s,\tilde s\in [0,\tau), s<\tilde s,$ such that $\tau=k(\tilde s-s)$ and $\tilde \tau =\tilde s-s$ is the smallest real number $\tau$ such that (\ref{zhaoh61}) holds. 

\item{\underline {Case (iii)}}. For any $s, t\in \mathbb R$, $\rho_s=\rho_t$. So $\bar \rho=\rho_s$ is an invariant measure for $\{P(t)\}_{t\geq 0}$. 
\vskip5pt

\underline {Case (i)} implies the following equivalent statements:
\vskip5pt

(ia). There exists a nontrivial angle variable  with $\lambda ={2l\pi\over \tau}$ for some $l\in {\mathbb N}\setminus \{0\}$ and no nontrivial angle variables with $\lambda <{2l\pi\over \tau}$;
\vskip5pt

(ib). The infinitesimal generator ${\mathcal A}$ of $U_t$ has infinite many simple eigenvalues $\{{2lm\pi\over \tau}\}_{m\in {\mathbb  Z}}$ for some $l\in {\mathbb N}\setminus \{0\}$, and no other eigenvalues;
\vskip5pt
 
(ic). The infinitesimal generator ${\mathcal L}$ of the semigroup $P(t)$ has infinite many simple eigenvalues 
$\{{2lm\pi\over \tau}i\}_{m\in {\mathbb  Z}}$ for some $l\in {\mathbb N}\setminus \{0\}$, and no other eigenvalues
on the imaginary axis.
\\

\underline {Case (ii)} implies the following equivalent statements:
\vskip5pt

(iia). There exists a nontrivial angle variable  with $\lambda ={2l\pi\over {\tilde\tau}}$ for some $l\in {\mathbb N}\setminus \{0\}$ 
and no nontrivial angle variables with $\lambda <{2l\pi\over \tilde \tau}$ for some $\tilde \tau ={\tau\over k}$,  
with some $k\in {\mathbb N}\setminus \{0,1\}$;
\vskip5pt

(iib). The infinitesimal generator ${\mathcal A}$ of $U_t$ has infinite many simple eigenvalues 
$\{{2lm\pi\over \tilde \tau}\}_{m\in {\mathbb  Z}}$ for some $l\in {\mathbb N}\setminus \{0\}$ and  some $\tilde \tau ={\tau\over k}$,  
with some $k\in {\mathbb N}\setminus \{0,1\}$, and no other eigenvalues;
\vskip5pt

(iic). The infinitesimal generator ${\mathcal L}$ of the semigroup $P(t)$ has infinite many simple eigenvalues 
$\{{2lm\pi\over \tilde \tau}i\}_{m\in {\mathbb  Z}}$ for some $l\in {\mathbb N}\setminus \{0\}$ and  some $\tilde \tau ={\tau\over k}$,  
with some $k\in {\mathbb N}\setminus \{0,1\}$, and no other eigenvalues
on the imaginary axis. 
\\

\underline {Case (iii)} is equivalent to the following equivalent statements:
\vskip5pt

(iiia). The angle variable is a constant and $\lambda =0$;
\vskip5pt

(iiib). The infinitesimal generator ${\mathcal A}$ of $U_t$ has one simple eigenvalue $0$ and no other eigenvalues;
\vskip5pt

(iiic). The infinitesimal generator ${\mathcal L}$ of the semigroup $P(t)$ has only one simple eigenvalues $0$, and no other eigenvalues
on the imaginary axis. 
\vskip5pt

(iiid). There exists $s\in [0,\tau)$ and a sequence $s_k\to s$, $s_k>s$ such that $L_s\cap L_{s_k}\neq \emptyset$.

\vskip5pt

Conversely, 
if there exists a nontrivial angle variable  with $\lambda ={2\pi\over \tau}$
and no nontrivial angle variables with $\lambda <{2\pi\over \tau}$, then $\tau$ is the minimum period of the 
periodic measure; if there exists a nontrivial angle variable  with $\lambda ={2\pi\over {\tilde\tau}}$ 
and no nontrivial angle variables with $\lambda <{2\pi\over \tilde \tau}$ for some $\tilde \tau={\tau\over k}, k\in {\mathbb N}\setminus \{0\}$, the minimum period of the periodic measure is no less than $\tilde \tau$; if the angle variable is a constant and $\lambda =0$, then the periodic measure has no positive minimum period, i.e. the periodic measure is a stationary measure.

\vskip5pt
\end{thm}
\begin{proof}
 It is obvious 
that there are only 3 possible cases (i), (ii) and (iii). 
First assume that for each $s\in {\mathbb R}$, $\rho_s$ as an invariant measure of $\{P(k\tau)\}_{k\in {\mathbb N}}$ is mixing.

\vskip5pt

\noindent
\underline {Case (i)}. Now we prove that (i) implies (ia). 

 First suppose (i) holds. Note as a special case of Theorem 3.4.2 in \cite{da-prato}, for any $\Gamma \in {\mathcal B}(\mathbb{X})$, 
 $$
 P(t+k\tau,x,\Gamma)\to \rho_t(\Gamma)\ \ {\rm in }\ L^2(L_0,\rho_0),
 $$
  as $k\to \infty$. But all the measures $\rho_t$ are different for different $t\in [0,\tau)$. It follows from applying Theorem 3.4.1 in \cite{da-prato} that the invariant measure $\bar\rho$ is definitely not weakly mixing.
 Thus by Koopman-von Neumann theorem, there is an 
 angle variable that is not constant. 
 Then by Remark \ref{zhao20157d}, there is an angle variable $\alpha$ such that (\ref{zhao20157b}) holds and 
 $\lambda \ne 0$ and (\ref{zhao20157c}) is satisfied. 
 By Proposition 3.2.1 in \cite{da-prato}, there exists a function $\phi\in L^2_{\mathbb C}(\mathbb{X}, d\bar\rho)$ such that 
 \begin{eqnarray*}
 P(t)\phi ={\rm e}^{i\lambda t}\phi, {\rm \ for\ any}\ t\geq 0,\ \bar\rho-a.s.,
 \end{eqnarray*}
 and $\xi$ defined in (\ref{zhao20157c}) is given by $\xi(\omega^*)=\phi(\omega^*(0))$.  In particular, there exists $s\in [0,\tau)$
 such that $\int_{L_s}(\phi(x))^2\rho_s(dx)>0$ and  
 \begin{eqnarray*}
 P(k\tau)\phi (x) ={\rm e}^{ik\lambda \tau}\phi (x),{\rm\ for\ any}\ k\in {\mathbb N},  \ x\in L_s.
 \end{eqnarray*}
 However, the discrete random dynamical system $\Phi(k\tau)$, by Remark \ref{zhao100m} (ii), 
 starting from $L_s$ will return on $L_s$ with probability $1$. Furthermore on $L_s$,  the invariant measure $\rho_s$ of $\Phi(k\tau)|_{L_s}$ is mixing. By Theorem 
3.4.1 in \cite{da-prato}, ${\rm e}^{ik\lambda \tau}=1$ and $\phi |_{L_s}$ is constant.
This suggests that $\lambda k\tau$ is divisible by $2\pi$ for any $k$. In particular, $\lambda \tau$ is divisible by $2\pi$, so $\lambda \tau =2l\pi$ for certain $l\in {\mathbb N}\setminus \{0\}$. We can certainly choose the smallest 
such $\lambda$, still denoted by $\lambda$ without causing any confusions.  
The claim that (i) implies (ia) is asserted.
  
 The equivalence of (ia) and (ib) follows from Remark \ref {zhao20157d}.

We now prove the equivalence of (ib) and (ic). If (ib) is true, then $U_t$ has eigenvalues ${\rm e}^{i{2m\pi\over \tau}t}, \ m\in {\mathbb Z}$. Thus by
the result that the eigenvalues of $P(t)$ on the unit circle are the same as the eigenvalues of $U_t$, so ${\rm e}^{i{2m\pi\over \tau}t}, \ m\in {\mathbb Z}$,
are only eigenvalues of $P(t)$ on the unit circle, and they are simple. Then it follows from the definition (\ref{zhao300a}) of ${\mathcal L}$, 
$\{i{2m\pi\over \tau}\}_{m\in {\mathbb Z}}$ are only simple eigenvalues of ${\mathcal L}$ on the imaginary axis. The converse can be proved similarly.

\vskip5pt

 \noindent
\underline {Case (ii)}. 
The proof that (ii) implies (iia)  and equivalence of (iia), (iib) and (iic)) can be done by a similar argument as in the proof in case (i).
 \vskip5pt
 
\noindent
\underline{Case (iii)}. 
The equivalence of (iii) and (iiia). The part from (iii) to (iiia) was already 
given when we consider the Case (i). Now we assume (iiia) holds. In this case $\bar \rho$ is weakly mixing. 
In both Case (i) and Case (ii), $\bar\rho$ is not weakly mixing. So Case (iii) must occur and (iii) holds.

The equivalence of (iiia) and (iiib) follows from Koopman-von Neumann theorem and the equivalence of (iiia) with 
$\bar \rho$ being weakly mixing. The proof of the equivalence of (iiib) and (iiic) can be done similarly as the proof of the equivalence of 
(ib) and (ic).

We finally prove that (iii) and (iiid) are equivalent.  
Suppose (iiid) is true, we need to prove that $\rho_t=\rho_s$ for any $s, t\in \mathbb R$. First note under the stochastic continuity assumption, it is well known that for any $\phi\in C_b(\mathbb{X})$,
\begin{eqnarray}\label{eq2.46}
\lim_{t\to 0} \int _\mathbb{X}\phi(y)P(t,x,dy)=\phi (x).
\end{eqnarray}
For each $k$, set $\tau_k=s_k-s$. Then $\tau_k\to 0$ as $k\to\infty$, and by 
Theorem \ref{zhao311s}, $\rho_{s+\tau_k}=\rho_s$. Define for any $t>s$, there exists 
$N_k\in \mathbb N$ and $0\leq \lambda_k<\tau_k$ such that $t=s+N_k\tau_k+\lambda _k$. 
It is obvious that $\lambda_k\to 0$ as $k\to \infty$. So by the Chapman-Kolmogorov equation, (\ref{eq2.46}) and Lebesgue's dominated convergence theorem,
\begin{eqnarray*}
&&<\phi,\rho_t>\\
&=& \int _\mathbb{X}\phi(y)\rho_t(dy)\\
&=& \int _\mathbb{X}\phi(y)\rho_{s+N_k\tau_k+\lambda _k}(dy)\\
&=& \int _\mathbb{X}\phi(y)\int_\mathbb{X}P(\lambda _k,x,dy)\rho_{s+N_k\tau_k}(dx)\\
&=&  \int _\mathbb{X}\phi(y)\int_\mathbb{X}P(\lambda_k,x,dy)\rho_{s}(dx)\\
&=&  \int _\mathbb{X}(\int_\mathbb{X}\phi(y)P(\lambda_k,x,dy))\rho_{s}(dx)\\
&\to&  \int_\mathbb{X} \phi(x)\rho_s (dx)\\
&=&<\phi, \rho_s>.
\end{eqnarray*}
So $<\phi,\rho_s>=<\phi,\rho_t>$ for any $\phi \in C_b(\mathbb{X})$. Thus $\rho_t=\rho_s$. The result (iii) is proved.

The converse part that (iii) implies (iiid) is trivial.
\vskip5pt

Now we prove the converse part of the theorem. We assume there exists a nontrivial angle variable  with $\lambda ={2\pi\over \tau}$
and no nontrivial angle variables with $\lambda <{2\pi\over \tau}$. Note from Remark \ref{zhao20157d}, (\ref{zhao20157b}) is always true since $\bar\rho$ is ergodic. We now prove (i) 
 by contradiction. If $\tau$ is not the smallest number such that (\ref{zhaoh61}) holds, then either Case (ii) or Case (iii) should happen. 
 If Case (ii) happens, then by a similar argument as in the last paragraph, we can show that $\lambda ={2l\pi\over \tilde \tau}$, 
 and no nontrivial angle variables with $\lambda <{2l\pi\over \tilde \tau}$, for certain $l\in {\mathbb N}\setminus \{0\}$, where $\tilde \tau$ is the number given in 
 (ii). 
 This is a contraction. If Case (iii) happens,   
 then $\bar\rho$ is equal to $\rho_s$ for any $s$ and is weakly mixing. The proof is completely 
 independent of any argument in this part, so we can use the result without causing any confusions.
 This then leads us to conclude that $\lambda =0$ following 
 the Koopman-von Neumann theorem. This is also a contradiction. 
 Claim (i) follows. 

Now assume there exists a nontrivial angle variable  with $\lambda ={2\pi\over {\tilde\tau}}$ 
and no nontrivial angle variables with $\lambda <{2\pi\over \tilde \tau}$ for some $\tilde \tau={\tau\over k}, k\in {\mathbb N}\setminus \{0\}$. If the minimum period of the periodic measure is less than $\tilde \tau$, say $\tilde \tau^*<\tilde \tau$. Then from the result in case (ii) that we have proved,
 there is an angle variable with  $\lambda ={2l\pi\over \tilde \tau^*}$ for some $l\geq 1$ and no any nontrivial angle variable with $\lambda <{ 2l\pi\over \tilde \tau^*}$. This is a contradiction with the assumption. The claim that the minimum period of the period measure 
 is no less than $\tilde \tau$ is proved. 
 
 The very last claim has been already proved in case (iii).
\end{proof}

Noting the relation of the eigenvalues of the infinitesimal generator ${\mathcal L}$ on the imaginary axis and the angle variable mentioned above already, we can state
the converse part of Theorem \ref{zhao20157f} differently.

\begin{cor}\label{zhao118}
 Assume the transition probability is stochastically continuous and 
has a periodic measure $\{\rho_s\}_{s\in {\mathbb R}}$ of period $\tau$, which is PS-mixing. 
If the infinitesimal generator ${\mathcal L}$ has simple eigenvalues  $\{{2m\pi\over \tau}i\}_{m\in {\mathbb Z}}$
and no other eigenvalues on the imaginary axis, then 
the period $\tau$ is the minimum period of the periodic measure; 
if the infinitesimal generator ${\mathcal L}$ has simple eigenvalues $\{{2m\pi\over \tilde \tau}i\}_{m\in {\mathbb Z}}$, where $\tilde \tau={\tau\over k}$ for some $k\in {\mathbb N}, k\geq 1$ and no other eigenvalues on the imaginary axis, then 
 the minimum period of the periodic measure is no less than $\tilde \tau$; if the infinitesimal generator ${\mathcal L}$ has simple eigenvalue $\{0\}$
and no other eigenvalues  
on the imaginary axis, 
  then the periodic measure has no positive minimum period, i.e. the periodic measure is a stationary measure.
\end{cor}

We can also present Theorem \ref{zhao20157f} as a sufficient-necessary condition to distinguish random periodic and 
stationary regimes.

\begin{thm} \label{zhao113}
 Assume the transition probability is stochastically continuous and 
has a periodic measure $\{\rho_s\}_{s\in {\mathbb R}}$ of period $\tau$, which is PS-mixing. 
Then the minimum period of the periodic measure is $\tilde \tau={\tau\over k}$, for certain $k\in {\mathbb N}\setminus \{0\}$,  if and only 
if that the infinitesimal generator ${\mathcal L}$ has simple eigenvalues $\{{2lm\pi\over  \tau}i\}_{m\in {\mathbb Z}}$, 
 for some $l\in {\mathbb N}\setminus \{0\}$, and no other eigenvalues  
on the imaginary axis.  The periodic measure has no positive minimum period if and only if that 
the infinitesimal generator ${\mathcal L}$ has simple eigenvalue $\{0\}$,
and no other eigenvalues  
on the imaginary axis.
\end{thm}
\begin{proof} 
The theorem follows from Theorem \ref{zhao20157f} and Corollary \ref{zhao118} easily.
\end{proof}

Dropping out the PS-mixing condition, the next theorem says that the PS-ergodicity can be obtained entirely based on the information of the spectral 
structure of the infinitesimal generator. Moreover, the Poincar\'e sections can also be defined by the eigenfunctions. This theorem improves significantly the result in the last part of Theorem \ref{zhao20157f}.

\begin{thm} \label{zhao114}
Assume the transition probability is stochastically continuous and 
has a periodic measure $\{\rho_s\}_{s\in {\mathbb R}}$ of period $\tau$.

(i). If the infinitesimal generator ${\mathcal L}$ has simple eigenvalues
$\{{2m\pi\over \tau}i\}_{m\in {\mathbb Z}}$,
and no other eigenvalues  
on the imaginary axis,  then the periodic measure 
$\{\rho_s\}_{s\in {\mathbb R}}$ is PS-ergodic and $\tau$ is the minimum period. Moreover, the eigenfunction 
$\phi_{0m}$ corresponding 
to the eigenvalue, $\lambda_m={2m\pi\over \tau}i, \ m=1,2,\cdots $,
is given by
\begin{eqnarray}\label{zhao1755}
\phi_{0m}(x)={\rm e}^{i{2m\pi\over \tau}t}, \ {when }\ x\in L_t.
\end{eqnarray}
Moreover, the Poincar\'e sections are given by the eigenfunction, 
denoted by $\phi_0$, corresponding to the eigenvalue ${2\pi\over \tau}i$, 
\begin{eqnarray}\label{zhao1754}
L_t=\{x\in {\mathbb X}: \phi_0(x)={\rm e}^{i{2\pi\over \tau}t}\},\ { for }\ t\in {\mathbb R}.
\end{eqnarray}

(ii). If the infinitesimal generator ${\mathcal L}$ has simple eigenvalues
$\{{2m\pi\over {\tilde \tau}}i\}_{m\in {\mathbb Z}}$, where $\tilde \tau={\tau\over l}$, for a $l\in {\mathbb N}$, 
and no other eigenvalues  
on the imaginary axis,  then the periodic measure 
$\{\rho_s\}_{s\in {\mathbb R}}$ is PS-ergodic and the minimum period of the invariant measure is at least $\tilde \tau$. 
\end{thm}

\begin{proof}
 (i). Let $\phi_0\in L^2_{\mathbb C}(L_0,\rho_0)$ satisfy 
\begin{eqnarray}\label{may112}
P(k\tau)\phi_0=\phi_0.
\end{eqnarray}
We will prove that $\phi_0$ is constant on $L_0$. Denote $\lambda =i{2\pi\over  \tau}$. Set  for $t\in {\mathbb R}$
\begin{eqnarray}
\phi_0^t(x):={\rm e}^{\lambda t}P(k\tau-t)\phi_0(x)={\rm e}^{\lambda t}\int_{L_0} P(k\tau-t,x,dy)\phi_0(y),\ \  x\in L_t,
\end{eqnarray}
where $k$ is the smallest integer such that $k\tau\geq t$. 
It is easy to know that 
\begin{eqnarray*}
\phi_0^{\tau}(x)={\rm e}^{\lambda \tau}\phi_0(x)=\phi_0(x),\ \ x\in L_0.
\end{eqnarray*}
Now by Jensen's inequality we see that 
$\phi_0^t\in L^2_{\mathbb C}(L_t,\rho_t)$ for each $t$.  It is easy to notice that $\{\phi_0^t\}_{t\in {\mathbb R}}$ is periodic in $t$. 
Moreover, it is noted that for any $s,t\geq 0$,
\begin{eqnarray}\label{may1}
P(s)\phi_0^{t+s}(x)&=&{\rm e}^{\lambda (t+s)}P(s)P(k\tau-(t+s))\phi_0(x)\nonumber\\
&=&{\rm e}^{\lambda (t+s)}P(k\tau-t)\phi_0(x)\nonumber\\
&=&{\rm e}^{\lambda s}{\rm e}^{\lambda t}P(k\tau-t)\phi_0(x)\nonumber
\\
&=&{\rm e}^{\lambda s}\phi_0^{t}(x),
\ \ x\in L_{t}.
\end{eqnarray}
Define
\begin{eqnarray*}
\phi_0(x)=\phi_0^t(x),\ \ {\rm for } \ x\in L_t, t\in {\mathbb R}.
\end{eqnarray*}
Then $\phi$ is well-defined on the whole space ${\mathbb X}$ and (\ref{may1}) is equivalent to
\begin{eqnarray}\label{zhao1752}
P(s)\phi_0={\rm e}^{\lambda s}\phi_0, \ {\rm for \ all} \ s\geq 0.
\end{eqnarray}
Thus
\begin{eqnarray}\label{may2}
{\mathcal L}\phi_0=\lambda \phi_0.
\end{eqnarray}
Now as the eigenvalue $\lambda $ of ${\mathcal L}$ is simple, so there is a unique, up to constant multiplication,
$\phi_0$ satisfying (\ref{may2}). However, it is observed that 
\begin{eqnarray}\label{zhao1756}
\phi_0(x)=\phi_0^t(x)={\rm e}^{\lambda t}, {\rm \ for }\ x\in L_t,
\end{eqnarray}
clearly satisfies (\ref{zhao1752}) and (\ref{may2}). In particular, $\phi_0(x)$ is constant on $L_0$.
Thus, 
$\rho_0$ is ergodic with respect to $\{P(k\tau)\}_{k\in {\mathbb N}}$. This means  the periodic measure is PS-ergodic.
Note that $\phi_0(x)$ are different when $x$ is in different Poincar\'e sections, and they are constant when $x$ is in a single 
Poincar\'e section. So $L_s\cap L_t=\emptyset $ when $s,t\in [0,\tau), s\ne t$. Thus $\tau$ is the minimum period. 
It is then obvious that $L_t$ can be constructed as (\ref{zhao1754}). 

Similarly, one can prove that the eigenfunction 
$\phi_{0m}$ corresponding 
to the eigenvalue, $\lambda_m={2m\pi\over \tau}i, \ m=1,2,\cdots $,
is given by (\ref{zhao1755}). 

(ii). Similar to the proof in (i), we also assume $\phi_0(x), x\in L_0$ 
satisfies (\ref{may112}). Consider $\lambda =i{2\pi\over  \tilde
\tau}$. Using the same procedure as above, one can construct the same eigenfunction as (\ref{zhao1756}), but with 
$\lambda$ given in this part. The eigenfunction also satisfies (\ref{zhao1752}) and (\ref{may2}). 
In particular, $\phi_0(x)$ is constant on $L_0$.
Thus, 
$\rho_0$ is ergodic with respect to $\{P(k\tau)\}_{k\in {\mathbb N}}$,
so the periodic measure is PS-ergodic.
Note that $\phi_0(x)$ are different when $x$ is in different Poincar\'e sections $L_t$ for $0\le t<\tilde \tau$, and they are constant when $x$ remains
 in a single 
Poincar\'e section. So $L_s\cap L_t=\emptyset $ when $s,t\in [0,\tilde \tau), s\ne t$. Thus the minimum period of the periodic measure is at least  $\tilde \tau$. 
\end{proof}

\begin{prop}\label{zhao117}
Assume the transition probability is stochastically continuous and 
has a periodic measure $\{\rho_s\}_{s\in {\mathbb R}}$ of period $\tau$, which is PS-ergodic.
Then there exist $k\in \mathbb N\setminus \{0\}$, $s,\tilde s\in [0,\tau], s<\tilde s,$ such that $\tau=k(\tilde s-s)$ and $\tilde \tau =\tilde s-s$ is the smallest real number $\tau$ such that (\ref{zhaoh61}) holds if and only if
there exist $s, \tilde s \in [0, \tau]$, $s<\tilde s$ such that $L_s\cap L_{\tilde s}\neq \emptyset$ 
and $L_s\cap L_r= \emptyset$ for any $r\in (s,\tilde s)$.
\end{prop}
\begin{proof}
Assume there exist $s, \tilde s \in [0, \tau]$, $s<\tilde s$ such that $L_s\cap L_{\tilde s}\neq \emptyset$ and $L_s\cap L_r= \emptyset$ for any $r\in (s,\tilde s)$. Then by Theorem \ref{zhao311s}, we have $\rho_{\tilde s}=\rho_s$. Thus
\begin{eqnarray}\label{eq2.44}
\rho_{s+\tilde \tau}=\rho_{\tilde s}=\rho_s=\rho_{s+\tau},
\end{eqnarray}
where $\tilde s=s+\tilde \tau$. 
Now for any $\Gamma\in \mathcal B(\mathbb{X})$, by the Chapman-Kolmogorov equation and (\ref{eq2.44}),
\begin{eqnarray*}
\rho_s(\Gamma)&=&\rho_{s+\tau}(\Gamma)\\
&=& \int_\mathbb{X} P(\tau-\tilde \tau, x, \Gamma)\rho_{s+\tilde \tau} (dx)\\
&=& \int_\mathbb{X} P(\tau-\tilde \tau, x, \Gamma)\rho_{s} (dx)\\
&=& \rho_{s+\tau-\tilde \tau}(\Gamma)\\
&&\cdots\\
&=&\rho_{s+\tau-k\tilde \tau}(\Gamma),
\end{eqnarray*}
where $k>0$ is an integer (unique) such that $0\leq \tau-k\tilde \tau<\tilde \tau$. Thus $\rho_s=\rho_{s+\tau-k\tilde \tau}$.
Note if $\tau>k\tilde \tau$, then $s< s+\tau-k\tilde \tau<s+\tilde \tau=\tilde s$. Because $L_s={\rm supp}
(\rho_s)={\rm supp}(\rho_{s+\tau-k\tilde\tau})=L_{s+\tau-k\tilde\tau}$, so it contradicts with the assumption that $L_s\cap L_r= \emptyset$ for any $r\in (s,\tilde s)$. 
Thus by the contradiction argument, we conclude that $\tau=k\tilde \tau$. Note for any $s'\geq s$, $\Gamma \in \mathcal B(\mathbb{X})$,
\begin{eqnarray}
\rho_{s'+\tilde \tau}(\Gamma\nonumber 
)&=&\int_\mathbb{X} P(s'-s,x,\Gamma)\rho_{s+\tilde \tau}(dx)\nonumber\\
&=&\int_\mathbb{X} P(s'-s,x,\Gamma)\rho_{s}(dx)\nonumber\\
&=&\rho_{s'}(\Gamma).\label{eq2.45}
\end{eqnarray}
We now claim that $\tilde \tau>0$ is the smallest number such that (\ref{eq2.45}) holds. If this is not true, there exists $\tau'\in (0,\tilde \tau)$ such that for any $\Gamma\in \mathcal B(\mathbb{X})$,
$$\rho_{s'+\tau'}(\Gamma)=\rho_{s'}(\Gamma).$$ Let $m$ be an integer number such that $s+m\tilde \tau\geq s'$. So by the same argument as (\ref {eq2.45}), we know 
$$\rho_{s+m\tilde\tau+\tau'}=\rho_{s+m\tilde\tau} =\rho_s.$$
But 
\begin{eqnarray*}
\rho_{s+m\tilde\tau+\tau'}(\Gamma)
&=&\int_\mathbb{X} P(\tau',x, \Gamma)\rho_{s+m\tilde\tau}(dx)\\
&=&\int_\mathbb{X} P(\tau',x, \Gamma)\rho_{s}(dx)\\
&=&\rho_{s+\tau'}(\Gamma).
\end{eqnarray*}
Thus 
$$\rho_{s+\tau'}=\rho_s.$$
This again is in contradiction with $\rho_r\neq \rho_s$ when $s<r<s+\tilde\tau$. Similar as above we can prove that when $s'<s$, (\ref{eq2.45}) is also true and $\tilde\tau$ is the smallest number for such 
an equality for all $s\in \mathbb R$. 

Conversely, if there exist $k\in \mathbb N\setminus \{0\}$, $s,\tilde s\in [0,\tau], s<\tilde s,$ such that $\tau=k(\tilde s-s)$ and $\tilde \tau =\tilde s-s$ is the smallest real number $\tau$ such that (\ref{zhaoh61}) holds,  it is trivial that there exist $\tilde s,s\in [0,\tau], \tilde s>s$ such that $\tilde \tau =\tilde s-s$ and $L_s=L_{\tilde s}$. The result then follows from Theorem
\ref {zhao311s}. 
\end{proof}

\section
{Spectral gap and PS-mixing}

%
%
%

We further this study here to prove that a spectral gap of the semigroup implies the convergence of the 
transition probability to the periodic measure along the subsequence $\{P(k\tau)\}_{k\in {\mathbb N}}$ on Poincar\'e sections.
Under the spectral gap assumption,
we obtain that the periodic measure $\{\rho _s\}_{s\in {\mathbb R}}$ is PS-mixing and the mixing rate. 
 
Assume $e_m\in L^2(\mathbb{X},\bar\rho(dx))$ is the eigenfunction of ${\mathcal L}$ with $||e_m||_{L^2(\mathbb{X},\bar\rho(dx))}\linebreak
=1$ corresponding to the eigenvalue $\lambda _m={2m\pi\over \tau}i$ on the imaginary axis, for each $m\in {\mathbb Z}$. It is well-known that $e_0=1$. Define
\begin{eqnarray*}
\hat {\mathbb H}={\rm span} \{e_m, m\in {\mathbb Z}\}\subset L^2(\mathbb{X},\bar\rho(dx)),
\end{eqnarray*} 
and 
\begin{eqnarray*}
\tilde {\mathbb H}=\hat {\mathbb H}^{\perp}=\{f\in L^2(\mathbb{X},\bar \rho(dx)), <f,e_m>=0, m\in {\mathbb Z}\},
\end{eqnarray*} 
where $<f,g>=<f,g>_{L^2}$.

Consider 
$$
Q_s(k\tau):=P(k\tau)|_{L_s}: L^2(L_s,\rho_s(dx))\to L^2(L_s,\rho_s(dx)).  
$$ 
We say the discrete semigroup $\{P(k\tau)\}_{k\in {\mathbb N}}$, has spectral gap or is of exponential contraction on the Poincar\'e section $L_s$ if there exists a $\delta>0$ such that 
\begin{eqnarray}
\lim_{k\to \infty}{1\over k\tau}\ln ||Q_s(k\tau)|_{\tilde {\mathbb H}}||_s< -\delta <0, 
\end{eqnarray}
where $||Q_s(k\tau)||_s$ is the operator norm of $Q_s(k\tau)$
on $ L(L^2(L_s,\rho_s(dx))\cap \tilde {\mathbb H})$.

We prove the following result.

\begin{prop}\label{zhao300q}
Assume the Markovian semigroup $\{P(t)\}_{t\geq 0}$, has a periodic measure $\{\rho_s\}_{s\in {\mathbb R}}$ of period $\tau>0$, 
and the corresponding 
infinitesimal operator $\mathcal {L}$ has simple eigenvalues $\lambda _m={2m\pi\over \tau}i, m\in {\mathbb Z}$ only on the imaginary axis. 
If 
the semigroup $\{P(k\tau)\}_{k\in {\mathbb N}},$ has a spectral gap on the Poincar\'e 
section $L_0$, then for each $s\in \mathbb R$, the invariant measure $\rho_s$ of $\{P(k\tau)\}_{k\in {\mathbb N}}$
is mixing and for any $f\in L^2(\mathbb{X},\bar\rho(dx))$, we have that for a.e. $s\geq 0, k\in {\mathbb N}$, 
\begin{eqnarray}\label{zhao300m}
\ \ \ \ \ \ ||P(k\tau +s)f- \int _{L_s}f(x)\rho_s(dx)||_{L^2(L_0,\rho_0(dx))}
\leq {\rm e}^{-\delta k\tau}||f||_{L^2(L_s,\rho_s(dx))}.
\end{eqnarray}
Moreover, if the semigroup $\{P(k\tau)\}_{k\in {\mathbb N}},$ 
has a spectral gap on each Poincar\'e 
section $L_t$ for $t\in [0,\tau)$, then the periodic measure $\{\rho_s\}_{s\in {\mathbb R}}$ is PS-mixing, has minimum period $\tau$ 
and for any $f\in L^2(\mathbb{X},\bar\rho(dx))$, $k\in {\mathbb N}$,
\begin{eqnarray}\label{zhao300s}
 &&  \int _{\mathbb{X}}|{1\over \tau}\int _{k\tau}^{(k+1)\tau}P(t)f(x)dt- \int _{\mathbb{X}}f(x)\bar \rho(dx)|\bar\rho(dx)
\leq {\rm e}^{-\delta k\tau}||f||_{L^2(\mathbb{X},\bar\rho(dx))}.
\end{eqnarray}

\end{prop}
\begin{proof}
By 
the spectral gap assumption of the semigroup $P(t)$ on the Poincar\'e section $L_0$, it is easy to see that   $\rho_0$ as
the invariant measure of $\{P(k\tau)\}_{k\in {\mathbb N}}$ on $L_0$ is mixing, and 
for any $f\in \tilde {\mathbb H}\cap L^2(L_0,\rho_0(dx))$
\begin{eqnarray}\label{zhao300j}
||P(k\tau)f||_{L^2(L_0,\rho_0(dx))}\leq {\rm e}^{-\delta k\tau} ||f||_{L^2(L_0,\rho_0(dx))}.
\end{eqnarray}

For any  $f\in L^2(\mathbb{X},\bar\rho(dx))$, it is easy to see that $f\in L^2(L_s,\rho_s(dx))$ for a.e. $s\in [0,\tau)$. Consider for any fixed $t, s\in [0,\tau)$ and $f\in L^2(L_{t+s},\rho_{t+s}(dx))$, note by Jensen's inequality
\begin{eqnarray}
||P(s)f||_{L^2(L_t,\rho_t(dx))}
&=& \left [\int _{L_t}(\int _\mathbb{X} P(s,x,dy)f(y))^2\rho_t(dx)\right ]^{\frac 12}\nonumber\\
&\leq &\left [\int _{L_t}\int _\mathbb{X} P(s,x,dy)f^2(y)\rho_t(dx)\right ]^{\frac 12}\nonumber\\
&=&\left [\int _{L_{t+s}}f^2(y)\rho_{t+s}(dy)\right ]^{\frac 12}\nonumber\\
&=&||f||_{L^2(L_{t+s},\rho_{t+s}(dx))}.\label{eqn1.67}
\end{eqnarray}
so $P(s)f\in L^2(L_t,\rho_t(dx))$ and 
 there exist $\widehat {P(s)f}\in \hat {\mathbb H}$ and $\widetilde {P(s)f}\in \tilde {\mathbb H}$ such that 
$$
P(s)f=\widehat {P(s)f}+\widetilde {P(s)f}.
$$
Here 
$$
\widehat {P(s)f}=\sum _{m\in {\mathbb Z}}<e_m,P(s)f>e_m.
$$
By (\ref{zhao300j}), we derive that for any $f\in L^2(L_s,\rho_s(dx))$,
\begin{eqnarray}\label{zhao300k}
||P(k\tau)[P(s)f-\widehat {P(s)f}]||_{L^2(L_0,\rho_0(dx))}\leq 
{\rm e}^{-\delta k\tau}||P(s)f-\widehat {P(s)f}||_{L^2(L_0,\rho_0(dx))}.
\end{eqnarray}
Note that 
for any $ k\in {\mathbb N}, s\geq 0$,
\begin{eqnarray}\label{zhao300l}
P(k\tau)\widehat {P(s)f}
= \sum _{m\in {\mathbb Z}}<e_m,P(s)f>{\rm e}^{{2m\pi\over \tau}k\tau i} e_m
= \sum _{m\in {\mathbb Z}}<e_m,P(s)f>e_m
= \widehat {P(s)f}.
\end{eqnarray}
That is to say that $\widehat {P(s)f}$ is an eigenfunction of $P(k\tau)$ corresponding to eigenvalue $1$.
 By Theorem \ref{zhao20157f}, $\rho_0$ is an ergodic invariant measure with respect to $\{P(k\tau)\}_{k\in {\mathbb N}}$, 
on $L_0$, so $\widehat {P(s)f}$ is constant on $L_0$ by
Theorem 3.2.4 in \cite{da-prato}.

Moreover,  from (\ref{zhao300k}) and (\ref{zhao300l}), we have 
\begin{eqnarray}\label{zhao300p}
||P(k\tau +s)f-\widehat {P(s)f}||_{L^2(L_0,\rho_0(dx))}
&\leq& {\rm e}^{-\delta k\tau}||P(s)f-\widehat {P(s)f}||_{L^2(L_0,\rho_0(dx))}\nonumber\\
&\leq& {\rm e}^{-\delta k\tau}||P(s)f||_{L^2(L_0,\rho_0(dx))}\nonumber\\
&\leq& {\rm e}^{-\delta k\tau}||f||_{L^2(L_s,\rho_s(dx))}.
\end{eqnarray}
Now note that  $\widehat {P(s)f}$ is constant on $L_0$,
so by Jensen's inequality and (\ref{zhao300p}), we have
$$\int _{L_0}(P(k\tau +s)f)(x)\rho_0(dx)\to \widehat {P(s)f},
$$
as $k\to\infty$. 
However, by Fubini theorem and (\ref{zhaoh61}), 
$$
\int _{L_0}(P(k\tau +s)f)(x)\rho_0(dx)=\int _{L_s}f(x)\rho_s(dx).
$$ 
Thus $\widehat {P(s)f}=\int _{L_s}f(x)\rho_s(dx)$ and so (\ref{zhao300m}) holds for any $f\in L^2(L_s,\rho_s(dx))$.


If the semigroup $\{P(k\tau)\}_{k\in {\mathbb N}}$ has spectral gap on each Pincar\'e section, it is easy to see that 
the periodic measure $\{\rho_s\}_{s\in {\mathbb R}}$ is PS-mixing. 
Similarly, (\ref{zhao300m}) holds for $f\in L^2(L_{t+s},\rho_{t+s}(dx))$ i.e.
\begin{eqnarray}\label{zhao300t}
\int _{L_t} [P(k\tau +s)f(x)- \int _{L_{s+t}}f(x)\rho_{t+s}(dx)]^2\rho_t(dx)
\leq {\rm e}^{-2\delta k\tau}||f||^2_{L^2(L_{t+s},\rho_{t+s}(dx))}.
\end{eqnarray}
But for any $f\in L^2(\mathbb{X},\bar\rho(dx))$, we know $f\in L^2(\mathbb{X},\rho_t(dx))$ for a.e. $t$. In particular we have (\ref{zhao300t}) 
for a.e. $t\in [0,\tau)$. In particular 
$\rho_t$ is mixing with respect to $\{P(k\tau)\}_{k\in {\mathbb N}}$ on $L_t$. 
So it follows from applying Fubin's theorem, Jensen's inequality and (\ref{zhao300t}) that 
\begin{eqnarray*}
&&
\int _{\mathbb{X}}\left |{1\over \tau}\int _{k\tau}^{(k+1)\tau}P(s)f(x)ds- \int _{\mathbb{X}}f(x)\bar \rho(dx)\right |\bar\rho(dx)\\
&=&
\int _{\mathbb{X}}\left |{1\over \tau}\int _{0}^{\tau}P(k\tau+s)f(x)ds- {1\over \tau}\int _0^{\tau}\int _{\mathbb{X}}f(x)\rho_s(dx)ds\right |\bar\rho(dx)\\
&\leq&
{1\over \tau^2}
\int _0^{\tau}\int _{\mathbb{X}}\int _{0}^{\tau}\left |P(k\tau+s)f(x)-\int _{\mathbb{X}}f(x)\rho_{s+t}(dx)\right |ds\rho_t(dx)dt\\
&=&
{1\over \tau^2}
\int _0^{\tau}\int _{0}^{\tau}\int _{\mathbb{X}}\left |P(k\tau+s)f(x)-\int _{\mathbb{X}}f(x)\rho_{s+t}(dx)\right |\rho_t(dx)dsdt
\\
&=&
{1\over \tau^2}
\int _0^{\tau}\int _{0}^{\tau}\left [\int _{\mathbb{X}}\left |P(k\tau+s)f(x)-\int _{\mathbb{X}}f(x)\rho_{s+t}(dx)\right |^2\rho_t(dx)\right ]^{1\over 2}
dsdt
\\
&\leq &
{1\over \tau}
\int _0^{\tau}
{\rm e}^{-\delta k\tau}||f||_{L^2(L_t,\rho_t(dx))}dt
\\
&\leq &{\rm e}^{-\delta k\tau}||f||_{L^2(\mathbb{X},\bar\rho(dx))}.
\end{eqnarray*}
The proof is completed.
\end{proof}

\begin{thm} \label{zhao311m}
Assume the same conditions as in Proposition \ref{zhao300q}.
 Then the periodic measure is ergodic and for any $\Gamma \in {\mathcal B}(\mathbb{X})$, 
\begin{eqnarray}\label{zhao300r}
\int _{{\mathbb R}^d}|{1\over \tau}\int _{k\tau}^{(k+1)\tau}
P(s,x, \Gamma)ds- \bar\rho(\Gamma)|\bar\rho(dx)
\leq {\rm e}^{-\delta k\tau}.
\end{eqnarray}
\end{thm}

\begin{proof} The result follows from taking $f=I_{\Gamma}$, where $\Gamma \in {\mathcal B}({\mathbb R}^d)$ in (\ref{zhao300s}). 
Thus Condition A is satisfied with exponential convergence and so the periodic measure is ergodic from Lemma \ref{fengc4}.
\end{proof}

\section{Construction of random periodic paths from a periodic measure}

In general, with the original probability space, similar to the case that
an invariant measure does not give a stationary process, neither a periodic measure 
gives a random periodic path. In the following, an enlarged probability space and an extended 
random dynamical system will be constructed such that on the enlarged probability space,
a pull-back flow is a random periodic path of the extended random dynamical system. 
This construction is much more demanding than constructing
the periodic measure from a random periodic path. 

Now we consider a Markovian random dynamical system.  If it has a periodic measure on $(\mathbb{X},{\mathcal{B}}(\mathbb{X}))$, 
then we can construct a periodic measure on the product measurable space $(\Omega\times \mathbb{X}, {\mathcal{F}}\otimes {\mathcal{B}}(\mathbb{X}))$. 
Here we use Crauel's construction of invariant measures on the product space from invariant measures of transition 
semigroup on phase space (\cite{hans2}). 

\begin{thm}\label{zhao16} Assume the Markovian random dynamical system $\Phi$ has a periodic measure $\rho: \mathbb R\to {\mathcal
P}(\mathbb{X})$ on $(\mathbb{X},{\mathcal{B}}(\mathbb{X}))$. Then for any $s\in \mathbb{R}$
\begin{eqnarray}\label{zhaoh3}
(\mu _s)_{\omega}: =\lim_{n\to \infty}\Phi (n\tau+s,\theta({-n\tau-s})\omega)\rho_0,
\end{eqnarray}
exists. Let
$$
\mu_s(dx,d\omega)=(\mu_s)_{\omega}(dx)\times P(d\omega).
$$ Then $\mu_s$ is a periodic measure on the product measurable space $(\Omega\times \mathbb{X}, {\mathcal{F}}\otimes {\mathcal{B}}(\mathbb{X}))$ for $\Phi$ and $E(\mu_s)_{\cdot}=\rho_s$, $s\in \mathbb{R}$.
\end{thm}

\begin{proof} 
First note that $\rho_0$ is a forward invariant measure under $P^*(n\tau),\ n\in \mathbb{N}$. By Crauel \cite{{hans2}}, we know that the following limit exists
$$({\mu_0})_\omega:=\lim_{n\to \infty}\Phi (n\tau,\theta({-n\tau})\omega)\rho_0.$$
By cocycle property of $\Phi$, we have that for any $B\in {\mathcal{B}}(\mathbb{X})$ for any $s\in \mathbb{R}^+$,
\begin{eqnarray}
&&
\lim_{n\to \infty}\Phi (n\tau+s,\theta({-n\tau-s})\omega)\rho_0(B)\nonumber\\
&=&\lim_{n\to \infty}(\Phi(s,\theta({-s})\omega)\circ \Phi (n\tau,\theta({-n\tau})\theta({-s})\omega)\rho_0)(B)\nonumber\\
&=&\lim_{n\to \infty}(\Phi (n\tau,\theta({-n\tau})\theta({-s})\omega)\rho_0)(\Phi(s,\theta({-s})\omega)^{-1}B)\nonumber\\
&=&(\mu_0)_{\theta(-s)\omega}(\Phi(s,\theta({-s})\omega)^{-1}B)\nonumber\\
&=&\Phi(s,\theta({-s})\omega)(\mu_0)_{\theta(-s)\omega}(B)\nonumber\\
&=& :(\mu_s)_{\omega}(B).\label{15}
\end{eqnarray}
When $s\in \mathbb{R}^-$, we can also obtain that the above limit still exists by decomposing $s=-m\tau +s_0$, $s_0\in [0,\tau)$,  and considering
\begin{eqnarray*}
&&
\lim_{n\to \infty}\Phi (n\tau+s,\theta({-n\tau-s})\omega)\rho_0(B)
\nonumber
\\
&=&\lim_{n\to \infty}(\Phi(s+m\tau ,\theta({-(s+m\tau )})\omega)\\
&&\ \ \ \ \ \ \ \ \ \ \circ \Phi ((n-m)\tau,\theta(-(n-m)\tau)\theta (-(s+m\tau) )\omega)\rho_0)(B)
\nonumber\\
&=& :(\mu_s)_{\omega}(B).
\end{eqnarray*}
Now, from the cocycle property and (\ref{zhaoh3}) and the argument of taking limits in (\ref{15}), we know that for $t\in \mathbb{R}^+$,
\begin{eqnarray}\label{zhao15}
\Phi(t,\omega)(\mu_s)_\omega
&=&\lim_{n\to \infty}\Phi(t,\omega)\circ \Phi (n\tau+s,\theta({-n\tau-s})\omega)\rho_0\nonumber\\
&=&\lim_{n\to \infty} \Phi (n\tau+t+s,\theta({-n\tau-t-s})\theta(t)\omega)\rho_0\nonumber\\
&=&(\mu_{t+s})_{\theta(t)\omega}.
\end{eqnarray}
It then follows from a standard argument that for any $A\in {\mathcal{F}}\otimes{\mathcal{B}}(\mathbb{X})$, by (\ref{zhao12})
and (\ref{zhao15}), for $t\in \mathbb{R}^+$
\begin{eqnarray*}\label{1.5}
(\bar\Theta(t)\mu_s)(A)
=
\mu_{t+s}(A).
\end{eqnarray*}
Moreover, it is easy to see that
$$(\mu_{s+\tau})_\omega=\lim_{n\to \infty}\Phi ((n+1)\tau+s,\theta({-(n+1)\tau-s})\omega)\rho_0=(\mu_{s})_\omega,$$
so $$\mu_{s+\tau}=\mu_s.$$
Then $\mu_.$ is a periodic measure on the product measurable space $(\Omega\times \mathbb{X}, {\mathcal{F}}\otimes {\mathcal{B}}(\mathbb{X}))$ for $\Phi$.

Next let us prove for any $B\in {\mathcal{B}}(\mathbb{X})$,  $s\in \mathbb{R}$, $E(\mu_s)_\omega(B)=\rho_s(B)$. First, we will show that for any $B\in {\mathcal{B}}(\mathbb{X})$, $E(\mu_0)_\omega(B)=\rho_0(B)$.
In fact, by the Lebesgue's dominated convergence theorem, the Fubini theorem and measure preserving property of $\theta$,
\begin{eqnarray*}
E(\mu_0)_\omega(B)
&=&\int _{\Omega}\lim_{n\to \infty}\Phi (n\tau,\theta({-n\tau})\omega)\rho_0(B)P(d\omega)\\
&=&\lim_{n\to \infty}\int _{\Omega} \rho_0(\Phi (n\tau,\theta({-n\tau})\omega)^{-1}(B))P(d\omega)\\
&=&\lim_{n\to \infty}\int_\Omega\int_\mathbb{X} I_{\Phi (n\tau,\theta({-n\tau})\omega)^{-1}B}(x)d\rho_0(x)P(d\omega)\\
&=&\lim_{n\to \infty}\int_\mathbb{X} \int_\Omega I_{\Phi (n\tau,\theta({-n\tau})\omega)^{-1}B}(x)P(d\omega)d\rho_0(x)\\
&=&\lim_{n\to \infty}\int_\mathbb{X} \int_\Omega I_B (\Phi (n\tau,\theta({-n\tau})\omega)x)P(d\omega)d\rho_0(x)\\
&=&\lim_{n\to \infty}\int_\mathbb{X} P(n\tau, x, B)d\rho_0(x)\\
&=&\rho_0(B).
\end{eqnarray*}
Similarly and also applying the above result, we have for $s\in \mathbb{R}^+$,
\begin{eqnarray*}
E(\mu_s)_\omega(B)
&=&\lim_{n\to \infty}\int_\mathbb{X} P(s+n\tau, x, B)d\rho_0(x)\\
&=&\lim_{n\to \infty}\int_\mathbb{X}\int_\mathbb{X} P(s,y,B)P(n\tau, x, dy)\rho_0(dx)\\
&=&\lim_{n\to \infty}\int_\mathbb{X}P(s,y,B)\int_\mathbb{X} P(n\tau, x, dy)\rho_0(dx)\\
&=&\int_\mathbb{X} P(s,y,B)d\rho_0(dy)\\
&=&\rho_s(B).
\end{eqnarray*}
If $s\in \mathbb{R}^-$, there exists $m\in Z^+$, $s_0\in [0,\tau)$ such that $s=-m\tau+s_0$. So similarly as above
\begin{eqnarray*}
E(\mu_s)_\omega(B)
&=&\lim_{n\to \infty}\int_\Omega\Phi (s_0+(n-m)\tau,\theta({-s_0-(n-m)\tau})\omega)\rho_0(B)P(d\omega)\\
&=&\lim_{n\to \infty}\int_\mathbb{X} P(s_0+(n-m)\tau, x, B)d\rho_0(x)\\
&=&\rho_{s_0}(B)\\
&=&\rho_{s}(B).
\end{eqnarray*}
In summary, we proved the last claim of the theorem for all $s\in \mathbb{R}$.
\end{proof}

 We assume that the cocycle $\Phi$ generates a periodic probability measure $\mu$ on the product 
 measurable space $(\bar\Omega, \bar{\mathcal {F}})=(\Omega\times \mathbb X, {\mathcal{F}}\otimes{\mathcal{B}}(\mathbb{X}))$.
 The following observation of an extended probability space, a random dynamical system and the
 correct construction of an invariant measure  $\hat \mu$
are key to the proof of the following theorem, which enables us to construct  periodic paths from periodic measures.

Set $I_{\tau}=[0,\tau)$ of the additive modulo $\tau$, ${\mathcal{B}}(I_{\tau})=\{\emptyset, I_{\tau}\}$, $\hat\Omega=I_{\tau}\times \Omega\times \mathbb{X}$, $\hat{\mathcal{F}}={\mathcal{B}}
(I_{\tau})\otimes{\mathcal{F}}\otimes{\mathcal{B}}(\mathbb{X})$, $\hat\omega=(s,\omega,x)\in \hat\Omega$. Define the skew product $\hat \Theta : \mathbb{R}^+\times \hat \Omega\to \hat \Omega$ as
\begin{eqnarray}\label{zhaoh1}
\hat\Theta(t)\hat\omega
&=&(s+t\  mod \ \tau, \theta(t)\omega, \Phi(t,\omega)x)\nonumber\\
&=&
(s+t-[{{s+t}\over \tau}]\tau,\theta(t)\omega,\Phi(t,\omega)x),\ \  t\in \mathbb{R}^+.
\end{eqnarray}

\begin{thm}\label{thm1.3}
Assume that a random dynamical system $\Phi$ generates a periodic probability measure $\mu$ on the product measurable space $(\Omega\times \mathbb{X}, {\mathcal     {F}}\otimes {\mathcal     {B}}(\mathbb{X}))$.
Then a measure $\hat \mu$ on the measurable space $(\hat\Omega, \hat{\mathcal{F}})$ defined by,
\begin{eqnarray}\label{prop6g}
\hat \mu(I_{\tau}\times A)={1\over {\tau}}\int _0^\tau  \mu_s(A)ds,\ \ \hat \mu(\emptyset\times A)=0,
\end{eqnarray}
for any $ A\in {\mathcal{F}}\otimes{\mathcal{B}}(\mathbb{X})$, 
is a probability measure and $\hat\Theta(t):\hat \Omega\to \hat\Omega$  defined by (\ref{zhaoh1}) is measure $\hat \mu$-preserving, and
\begin{eqnarray}\label{zhaoh2}
\hat\Theta(t_1)\hat\Theta(t_2)=\hat\Theta(t_1+t_2), \ {\rm for \ any}\ t_1,t_2\in \mathbb{R}^+.
\end{eqnarray}
If we extend $\Phi$ to a map
over the metric dynamical system $(\hat\Omega, \hat{\mathcal{F}}, \hat \mu,(\hat\Theta(t))_{t\in \mathbb{R}^+})$ by
\begin{eqnarray}
\hat\Phi(t, \hat\omega)=\Phi(t,\omega),\ t\in \mathbb{R}^+,
\end{eqnarray}
 then $\hat\Phi$ is a RDS on $\mathbb{X}$ over $\hat\Theta$ and has a random periodic path $\hat Y: \mathbb{R}^+\times \hat \Omega\to \mathbb{X}$ constructed as follows:  for any $\hat\omega^*=(s,\omega^*,x^*(\omega^*))\in \hat \Omega$,
\begin{eqnarray}\label{prop6f}
\hat Y(t,\hat\omega^*):=\Phi(t+s, \theta(-s)\omega^*)x^*(\theta(-s)\omega^*),\ t\in \mathbb{R}^+.
\end{eqnarray}

\end{thm}
\begin{proof} It is easy to see that the proof of (\ref{zhaoh2}) is a matter of straightforward computations and $\hat \mu$ is a probability measure. To verify
$\hat\Theta(t)\hat \mu=\hat \mu$, for any $t\in \mathbb{R}^+$, first using (\ref{1.3}) and a similar argument as (\ref{im1}), we have that for any $t\in [0,\tau)$,
\begin{eqnarray*}
\hat\Theta(t)\hat \mu(I_{\tau}\times A)
&=&{1\over {\tau}}\int _0^\tau  \mu_s(\bar \Theta^{-1}(t)A)ds\\
&=&{1\over {\tau}}\int _0^\tau \mu_{s+t}(A)ds\\
&=&\hat \mu(I_{\tau}\times A).
\end{eqnarray*}
It is  trivial to note that $\hat\Theta(t)\hat \mu(\emptyset\times A)=\hat \mu(\emptyset\times A)$. So $\hat\Theta(t)$ is $\hat \mu$-preserving for $t\in [0,\tau)$. This can be easily generalised to any $t\in \mathbb{R}^+$ using the group property of $\hat\Theta$.
Moreover, it is trivial to see that $\hat \Phi$ is a cocycle on $\mathbb{X}$ over $\hat \Theta$.
Again,  the construction of $\hat Y$ given by (\ref{prop6f}) is key to the proof, from which the actual proof itself is quite straightforward.
In fact, for $\hat\omega=(s,\omega,x)$, we have
$\hat Y(t, \hat\omega)=\Phi(t+s, \theta(-s)\omega)x.$
Moreover, for any $r, t\in \mathbb{R}^+$, we have by the cocycle property
\begin{eqnarray}\label{zhao24}
\hat\Phi(r, \hat\Theta(t)\hat\omega)\hat Y(t,\hat\omega)
=\Phi(r, \theta(t)\omega)\Phi(t+s,\theta(-s)\omega)x
=\Phi(r+t+s, \theta(-s)\omega)x=\hat Y(r+t, \hat \omega).
\end{eqnarray}
Note that $\hat\Theta(\tau)\hat\omega=(s, \theta(\tau)\omega, \Phi(\tau,\omega)x)$, so we have
by the cocycle property
\begin{eqnarray}\label{zhao25}
\hat Y(t, \hat\Theta(\tau)\hat\omega)
&=&\Phi(t+s, \theta(\tau-s)\omega)\Phi(\tau,\theta(-s)\omega)x\nonumber\\
&=&\Phi(t+s+\tau,\theta(-s)\omega)x
=\hat Y(\tau+t,\hat\omega).
\end{eqnarray}
The proof is completed.
\end{proof}

\begin{rmk}
It is not clear how to extend the definition of $Y$ to ${\mathbb R^-}$ in general. However, if the cocycle $\Phi(t,\omega): \mathbb{X}\to \mathbb{X}$ is invertible for any $t\in \mathbb{R}^+$ and $\omega\in \Omega$, for instance in the case of SDEs in a finite dimensional space with some suitable conditions, it is obvious to extend
$Y$ to $\mathbb R^-$.
\end{rmk}

One implication of Theorems \ref{zhao16} and \ref{thm1.3} is that starting from a periodic measure $\rho_s\in {\mathcal{P}}(\mathbb{X}), \rho_{s+\tau}
=\rho_s, s\in \mathbb R$. one can construct
a (enlarged) probability space $(\hat\Omega, \hat{\mathcal{F}}, \hat \mu)$ and  an extended random dynamical system, with which the pull-back of the 
random dynamical system is a random periodic path. In the following we will prove that the 
transition probability of $\hat \Phi(t,\hat\omega)x$ is actually the same as $P(t,x,\cdot)$ and the
law of the random periodic solution $\hat Y$ is $\rho_s$, i.e.
$$\hat {\mathcal L}(\hat Y(s,\cdot ))=\rho_s,{\rm \ for \ any \ } s\in \mathbb{R}^+.
$$
We call $\hat Y$ a random periodic process as its law is periodic. 


In the following, by $\hat E$ we denote the expectation on 
$(\hat \Omega,\hat{\mathcal{F}},\hat \mu)$.

\begin{lem}\label{zhao26}
Assume $\rho_s$ is a periodic measure of a Markovian random dynamical system $\Phi$. 
Let the metric dynamical system $(\hat \Omega,\hat{\mathcal{F}},\hat \mu,(\hat \Theta (t))_{t\in \mathbb{R}^+})$, the extended random dynamical system $\hat \Phi$
and the random periodic process $\hat Y$ be defined in Theorem \ref{thm1.3}. Then for any $B\in {\mathcal{B}}(\mathbb{X})$
\begin{eqnarray*}
\hat \mu\{\hat\omega: \hat Y(t,\hat\omega)\in B\}=\rho_t(B),
\end{eqnarray*}
and
\begin{eqnarray*}
\hat P(t,y, B)=\hat \mu \{\hat\omega: \hat \Phi (t, \hat\omega)y\in B\}=P(t,y,B).
\end{eqnarray*}
Thus $\rho_{\cdot}$ is a periodic measure of $\hat \Phi$ as well.
\end{lem}
\begin{proof}
Note in (\ref {prop6g}), for any $A\in {\mathcal{F}}\otimes{\mathcal{B}}(\mathbb{X})$, by the periodicity of $\mu_s$ and measure preserving property of $\theta$,
\begin{eqnarray}
\hat \mu(I_{\tau}\times A)
&=&{1\over {\tau}}\int _0^\tau  \mu_s(A)ds\nonumber\\
&=&{1\over {\tau}}\int _0^\tau  \mu_{-s}(A)ds\nonumber\\
&=& {1\over {\tau}}\int _0^\tau \int_\Omega (\mu_{-s})_\omega(A_\omega)P(d\omega)ds\nonumber\\
&=& {1\over {\tau}}\int _0^\tau \int_\Omega (\mu_{-s})_{\theta(-s)\omega}(A_{\theta(-s)\omega})P(d\omega)ds.\label{eqn7.11}
\end{eqnarray}
From the proof of Theorem \ref{thm1.3}, we know that, for any $\hat \omega=(s,\omega, x)$,
$$
\hat Y(t, \hat\omega)=\Phi(t+s, \theta(-s)\omega)x,
$$
is a random periodic process on the probability space  $(\hat\Omega, \hat{\mathcal{F}}, \hat \mu, 
(\hat\Theta(t))_{t\in \mathbb{R}^+})$.
Then for any $t\in \mathbb{R}^+$ and $B\in {\mathcal{B}}(\mathbb{X})$, by (\ref{eqn7.11}), (\ref{zhao15}) and definition of $\hat Y$,
\begin{eqnarray*}
\hat \mu(\hat\omega: \hat Y(t, \hat\omega)\in B)
&=&\int_{\hat \Omega} I_B(\hat Y(t, \hat\omega))\hat \mu(d\hat\omega)\\
&=&{1\over {\tau}}\int _0^\tau \int_\Omega\int _\mathbb{X} I_B(\hat Y(t, \hat\omega))(\mu_{-s})_{\theta(-s)\omega}(dx)P(d\omega)ds\\
&=&{1\over {\tau}}\int _0^\tau \int_\Omega\int _\mathbb{X} I_B( \Phi(t+s, \theta(-s)\omega)x)(\mu_{-s})_{\theta(-s)\omega}(dx)P(d\omega)ds\\
&=&{1\over {\tau}}\int _0^\tau \int_\Omega\int _\mathbb{X} I_B(y) [\Phi(t+s, \theta(-s)\omega)(\mu_{-s})_{\theta(-s)\omega}](dy)P(d\omega)ds\\
&=&{1\over {\tau}}\int _0^\tau \int_\Omega\int _\mathbb{X} I_B(y)(\mu_{t})_{\theta(t)\omega}(dy)P(d\omega)ds\\
&=&E[(\mu_t)_{\theta(t)\cdot}(B)]=E[(\mu_t)_{\cdot}(B)]\\
&=&
\rho_t(B).
\end{eqnarray*}
 Now we consider
 $\hat\Phi(t, \hat\omega)=\Phi(t,\omega)$, the extended random dynamical system on $\mathbb{X}$ 
 over the probability space $(\hat \Omega, \hat {\mathcal{F}},\hat \mu)$. For any $y\in \mathbb{X}$, note again $\hat \omega=(s,\omega,x)$,
 \begin{eqnarray*}
 \hat P(t,y,B)
 &=&\hat \mu(\hat \omega: \hat \Phi(t, \hat\omega)y\in B)\\
 &=&\int _{\hat \Omega}[I_B(\hat\Phi(t,\hat\omega)y)]\hat \mu (d\hat \omega)\\
 &=&{1\over {\tau}}\int _0^\tau \int_\Omega\int _\mathbb{X} I_B(\hat\Phi(t,\hat\omega)y)(\mu_s)_{\omega}(dx)P(d\omega)ds\\
 &=&{1\over {\tau}}\int _0^\tau \int_\Omega\int _\mathbb{X} I_B(\Phi(t,\omega)y)(\mu_s)_{\omega}(dx)P(d\omega)ds\\
 &=&\int_\Omega I_B(\Phi(t,\omega)y)P(d\omega)\\
 &=&P(t,y,B).
 \end{eqnarray*}
 The last claim follows easily from the above two results already proved.

\end{proof}

{\bf Acknowledgements}. 
We would like to acknowledge the financial supports of a Royal Society Newton fund grant (ref. NA150344) and an EPSRC Established Career Fellowship
to HZ (ref. EP/S005293/1).


\begin{thebibliography}{99}

\setlength{\itemsep}{-1.2mm}


\bibitem{ar} L. Arnold, {\it Random Dynamical Systems}, Springer-Verlag Berlin (1998).

\bibitem{blw} P. W. Bates, K.N. Lu and B.X. Wang,  Attractors of non-autonomous stochastic lattice systems in weighted spaces,
{\it Physica D}, Vol. 289 (2014), 32-50.

\bibitem{BenziStochRes} R. Benzi, G. Parisi, A. Sutera
and A. Vulpiani, Stochastic resonance in climatic change, \emph{Tellus},
Vol 34 (1982), 10-16.

\bibitem{cce} A. P. Carverhill, M. J. Chappell and K. D. Elworthy, Characteristic exponents for stochastic flows, In: {\it Lecture Notes in Mathematics}, 
Vol. 1158, Springer, Berlin (1986), pp52-80.
 
\bibitem{chekroun} M. D. Chekroun, E. Simonnet and M. Ghil, Stochastic climate dynamics: random
attractors and time-dependent invariant measures, {\it Physica D}, Vol. 240 (2011), 1685-1700.

\bibitem{chen} M.F. Chen,
{\it Eigenvalues, Inequalities, and Ergodic Theory}, 
Probability and its Applications, Springer-
Verlag, 2005.

\bibitem{chojn} A. Chojnowska-Michalik, Periodic distribution for linear equations with general additive noise, {\it Bull. Pol. Acad.
Sci. Math.}, Vol. 38 (1990), 23-33.

\bibitem{hans2} H. Crauel, Markov measures for random dynamical systems. {\it Stochastics Stochastics Rep}, Vol. 37 (1991), no. 3, 153-173.

\bibitem{clrs} A.M. Cherubini, J.S.W. Lamb, M. Rasmussen and Y. Sato, A random dynamical systems perspective on stochastic resonance, {\it Nonlinearity}, Vol. 30 (2017), 2835-2853.

\bibitem{da-prato} G. Da Prato and J. Zabczyk, Ergodicity for infinite dimensional systems, {\it London Mathematical Society Lecture Note Series,} 229, Cambridge University Press,
1996.

\bibitem{doob} J. L. Doob, Asymptotic properties of Markoff transition probabilities, {\it Trans. Amer. Math. Soc.}, Vol. 63 (1948), 394-421.

\bibitem{nagel} K.-J. Engel and R. Nagel, {\rm One Parameter Semigroups for Linear Evolution Equations}, Graduate Texts in Mathematics 194, Springer, 1995. 
 
 \bibitem{flz1} C.R. Feng, Y. Liu and H.Z. Zhao, Numerical approximation of random periodic solutions of stochastic differential equations, {\it Zeitschrift fur angewandte Mathematik und Physik}, 
 Vol. 68 (2017), article 119, 1-32.
 
  \bibitem{flz3} C.R. Feng, Y. Liu and H.Z. Zhao,
 Numerical analysis of the weak schemes of random periodic
solutions of stochastic differential equations, in preparations,
 
 \bibitem{flz2} C.R. Feng, Y.J. Liu and H.Z. Zhao, ARMA model for random periodic processes, in preparations.
 
\bibitem{wu1} C.R. Feng, Y. Wu and H.Z. Zhao, Anticipating random periodic solutions--I.  SDEs with multiplicative linear noise, {\it J. Funct. Anal.}, Vol. 271 (2016), 365-417. 



\bibitem{fzz} C. Feng, H. Zhao and J. Zhong, Existence of geometric ergodic periodic measures of stochastic differential equations, 2019,  {\it arXiv:1904.08091}, submitted.

\bibitem{fzz2} C. Feng, H. Zhao and J. Zhong, Expected exit time for time-periodic stochastic differential equations
and applications to stochastic resonance, preprint. 


\bibitem{feng-zhao-zhou} C.R. Feng, H.  Z. Zhao and B. Zhou,
Pathwise random periodic solutions of stochastic differential equations,
{\it J. Differential Equations}, Vol. 251 (2011), 119-149.

\bibitem{feng-zhao} C.R. Feng and H.Z. Zhao,
Random periodic solutions of SPDEs via integral equations and Wiener-Sobolev compact embedding,
{\it J. Funct. Anal.}, Vol. 262 (2012), 4377-4422.

\bibitem{zhao2014} C.R. Feng and H.Z. Zhao, Random periodic processes, periodic measures and strong law of large numbers, Preprint, 2014, {\it arxiv.org/pdf/1408.1897v2.pdf}.
 
\bibitem{majda-2} B. Gershgorin, A. J. Majda, A test model for fluctuation-dissipation theorems with time-periodic statistics,
{\it Physica D}, Vol. 239 (2010), 1741-1757.

\bibitem{Ha} R. Z. Has'minskii, {\it Stochastic Stability of Differential Equations,} Springer, Second Edition, 2012.

\bibitem{hl} W. Huang and Z. Lian, Horseshoe and periodic orbits for quasi-periodic forced systems, {\it  	arXiv:1612.08394}, 2016.
 
 \bibitem{klunger} M. Klunger, Periodicity and Sharkovsky's Theorem for random dynamical systems, {\it Stochastics and Dynamics}, Vol. 1 (2001), 299-338.
 
\bibitem{lian} P. Lian and H.Z. Zhao, Pathwise properties of random mappings,
In: {\it New Trends in Stochastic Analysis and Related Topics--volume  in honour of Professor K.D. Elworthy},
edited by H.Z. Zhao and A. Truman, World Scientific, 2012, pp. 227-300.

\bibitem{poincare} H. Poincar\'e, Memoire sur les courbes definier par une equation differentiate.
J. Math. Pures Appli., Vol. 3 (1881), 375-442; J. Math. Pures Appli., Vol. 3 (1882), 251-296; J. Math. Pures Appli., Vol. 4 (1885), 167-244; J. Math. Pures Appli., Vol. 4 (1886), 151-217.

\bibitem{rockner1} N. Rezvani Majid and M. R\"ockner,
    The structure of entrance laws for time-inhomogeneous Ornstein-Uhlenbeck Processes with L\'evy Noise in Hilbert spaces, 
   arXiv:1507.06093.
   
\bibitem{rw} L. C. G. Rogers and D. Williams, {\it Diffusions, Markov processes and martingales, Vol. 2, It$\hat{\rm o}$ calculus}, Cambridge University Press, 2nd Edition, 2000.



\bibitem{scheutzow} M. Scheutzow, Periodic behaviour of stochastic Brusselator in the mean-field limit, {\it Probab. Th. Rel. Fields}, Vol. 72 (1986), 425-462.



\bibitem{wang} B.X. Wang, Existence, stability and bifurcation of random complete and
periodic solutions of stochastic parabolic equations, {\it Nonlinear Analysis}, Vol. 103 (2014), 9-25.

\bibitem{wangf} F.Y. Wang, {\it Functional Inequalities, Markov Semigroup and Spectral Theory},
 Chinese Sciences Press, Beijing, New York, 2005.
 
\bibitem{knobloch} B. Weiss and E. Knobloch, A stochastic return map for stochastic differential equations, {\it J. Stat. Phys.}, Vol. 58 (1990), 863-883.


\bibitem{zh-zheng} H.Z. Zhao and Z. H. Zheng,
Random periodic solutions of random dynamical systems,
{\it J. Differential Equations}, Vol. 246 (2009), 2020-2038.

\end{thebibliography}
\end{document}